\documentclass[11pt]{article}
\usepackage{amsmath}
\usepackage{amsfonts}
\usepackage{graphicx}
\usepackage{amssymb}
\usepackage{leftidx}
\usepackage{color}

\newtheorem{condition**}{A*}
\newtheorem{condition***}{C*}
\newtheorem{condition*}{C}

\newtheorem{proposition}{Proposition}[section]

\newtheorem{theorem}{Theorem}[section]

\newtheorem{remark}{Remark}[section]

\begin{document}

\title{Global Solutions of Stochastic Stackelberg Differential Games under
Convex Control Constraint}
\author{Liangquan Zhang$^{1}$\thanks{%
L. Zhang acknowledges the financial support partly by the National Nature
Science Foundation of China(Grant No. 11701040, 11871010 \&61871058) and the
Fundamental Research Funds for the Central Universities (No. 2019XD-A11).
E-mail: xiaoquan51011@163.com.}, Wei Zhang$^{1}$, \\
{\small 1. School of Science} \\
{\small Beijing University of Posts and Telecommunications} \\
{\small Beijing 100876, China} }
\maketitle

\begin{abstract}
This paper is concerned with a Stackelberg stochastic differential game,
where the systems are driven by stochastic differential equation (SDE for
short), in which the control enters the randomly disturbed coefficients
(drift and diffusion). The control region is postulated to be convex. By
making use of the first-order adjoint equation (backward stochastic
differential equation, BSDE for short), we are able to establish the
Pontryagin's maximum principle for the leader's global Stackelberg solution,
within adapted open-loop structure and closed-loop memoryless information
one, respectively, where the term global indicates that the leader's
domination over the entire game duration. Since the follower's adjoint
equation turns out to be a BSDE, the leader will be confronted with a
control problem where the state equation is a kind of \emph{fully} coupled
forward-backward stochastic differential equation (FBSDE for short).

As an application, we study a class of linear-quadratic (LQ for short)
Stackelberg games in which the control process is constrained in a closed
convex subset $\Gamma $ of full space $\mathbb{R}^{m}$. The state equations
are represented by a class of fully coupled FBSDEs with projection operators
on $\Gamma $. By means of monotonicity condition method, the existence and
uniqueness of such FBSDEs are obtained. When the control domain is full
space, we derive the resulting backward stochastic Riccati equations.
\end{abstract}



\section{Introduction}

\label{sect:1}

H. von Stackelberg \cite{S1934} first introduced a hierarchical solution for
markets with leaders and followers in 1934 in order to obtain optimal
strategies in competitive economics, which is now known as the Stackelberg
equilibrium. The Stackelberg game is also know as the leader-follower game,
whose economic background can be derived from some markets where certain
companies have advantages of domination over others. Stackelberg strategies
are rational and optimal solutions for both two players. A Stackelberg game
can be described briefly as follows: the leader first announces his/her
strategy at the beginning of the game. Based on the knowledge of the
leader's action, the follower, makes an instantaneous reaction by optimizing
his/her own performance index. Then, by anticipating the optimal response of
the follower, the leader will select an optimal action to optimize his/her
cost functional on the rational reaction curve of the follower. Therefore,
one player must make a decision after the other player's decision is made.
The leader's optimal action and the follower's rational response constitute
a Stackelberg solution.

Since its nice structure and background in economy, there have been a great
deal of substantial works along this research direction. For instance, the
leader-follower's feature can be applied in many fields, such as the
newsvendor/wholesaler problem (\O ksendal et al. \cite{OSU}), the optimal
reinsurance problem (Chen and Shen \cite{CS2018}), the operations management
and marketing channel problem (Li and Sethi \cite{LS2017}) and the
principal-agent/optimal contract problem (Cvitani\'{c} and Zhang \cite{CZSp}%
). Particularly, the feedback and adapted feedback information structures
have been employed in treating supply chain management, marketing channel
management, and economics problems in \cite{CS2012JOTA, CS2012, DJLS2000,
HPS2009, HPSG2007, KT2007} (for more information see reference therein).
Besides, the celebrated Pontryagin's maximum principle for stochastic
differential games within the framework of BSDE can be found in Wang and Yu
\cite{WY2010, WY2012} and Yu \cite{Yu2015}.

We state a few of them related our work. First, Castanon and Athans \cite%
{CA1979} considered an LQ stochastic dynamic Stackelberg strategies in the
early, and obtained a feedback Stackelberg solution for two-person nonzero
sum game. Yong \cite{Yong2002} studied an indefinite LQ leader-follower
stochastic differential game with random coefficients and control-dependent
diffusion. The related Riccati equations for the follower and for the leader
are derived sequentially to obtain the state feedback representation of the
open-loop Stackelberg equilibrium points, moreover the sufficient conditions
for their solvability with deterministic coefficients are given in the
special case. Ba\c{s}ar et al. \cite{BBS2010} introduced the notion of mixed
leadership in nonzero-sum differential games where one player could act as
both leader and follower, depending on the control variable. \O ksendal et
al. \cite{OSU} investigated a general stochastic Stackelberg differential
game with delayed information, established the maximum principle, and
applied it to continuous-time newsvendor problems. Bensoussan et al. \cite%
{BCS2015} introduced several global solution concepts in terms of the
players' information patterns, and derived the maximum principle by means of
FBSDEs for the leader's global Stackelberg solution under the adapted
open-loop and adapted closed-loop memoryless information structure (details
see below) with non-controlled diffusion term. Meanwhile they investigate
the LQ case where the weight matrices in the cost functionals are positive
definite and controls do not entre into the diffusion term of the state
equation. Mukaidani and Xu \cite{MX2015} considered the Stackelberg games
for linear stochastic systems driven by It\^{o} differential equations with
multiple followers. The Stackelberg strategies, obtained by using sets of
cross-coupled algebraic nonlinear matrix equations, are developed under two
different settings: the followers act either cooperatively to attain Pareto
optimality or non-cooperatively to arrive at a Nash equilibrium. Li and Yu
\cite{LY2015} provided the solvability of a coupled FBSDEs under a
multilevel self-similar domination-monotonicity structure, then it is
employed to characterize the unique equilibrium of an LQ generalized
Stackelberg stochastic differential game with hierarchy in a closed form.
Huang et al. \cite{HSW} studied a controlled linear-quadratic-Gaussian large
population system combining major leader, minor leaders and minor followers.
The Stackelberg-Nash-Cournot (SNC for short) approximate equilibrium is
obtained by means of the combination of a major-minor mean-field game and a
leader-follower Stackelberg game, besides the feedback form of the SNC
approximate equilibrium strategy is constructed through coupled Riccati
equations.

In this paper, we shall study the Stackelberg games under two stochastic
settings (taken from Benssousan et al. \cite{BCS2015}). The first one is
\textit{adapted open-loop} (AOL) which can be states in summary as follows:
given the strategy $u$ claimed by the leader at the beginning of the game,
the follower wants to minimize his cost functional $J_{2}(u;v)$ associated
with the leader's strategy $u$ on the whole duration of the game. His
optimal response $v^{\ast }$ will be an adapted process such that $%
J_{2}(u;v^{\ast }\left( u\right) )\leq J_{2}(u;v\left( u\right) ).$ The
leader makes an instantaneous reaction $u^{\ast }$ by optimizing his/her
performance index on the rational reaction curve of the follower, i.e., $%
J_{1}(u^{\ast };v^{\ast }\left( u^{\ast }\right) )\leq J_{1}(u;v^{\ast
}\left( u\right) ),$ anticipating the follower's optimal response $v^{\ast }$%
. The pair $(u^{\ast },v^{\ast })$ is called an AOL solution of the
Stackelberg game. The other one, \textit{adapted closed-loop memoryless}
(ACLM), comparing with AOL, the leader's strategy and the follower's
response strategy depend on the state (\textit{feedback form}), which turns
the control problem into a non-standard one. For leader's each strategy $u$
made in advanced, the follower would like to seek his optimal response $%
v^{\ast }$ such that $J_{2}(u;v^{\ast }\left( u\right) )\leq J_{2}(u;v\left(
u\right) ).$ Then, by taking the rational response of the follower into
account, the leader, of course, pick an action $u^{\ast }$ such that $%
J_{1}(u^{\ast };v^{\ast }\left( u^{\ast }\right) )\leq J_{1}(u;v^{\ast
}\left( u\right) ).$ The pair $(u^{\ast },v^{\ast })$ is called an ACLM
solution of the Stackelberg game.

To summarize the above, we see that the novelty of the formulation in this
paper is the following:

\begin{itemize}
\item Comparing with Bensoussan et al. \cite{BCS2015}, our diffusion term in
stochastic system allows to depend on control variable. As we shall see
Section \ref{sec2} below, due to this nice structure of our control system,
on the one hand, the adjoint equation for leader becomes more general no
matter of AOL or ACLM cases; on the other hand, the related stochastic
Riccati equation for the follower and the leader considered simultaneously
by putting the follower's Hamiltonian system as the leader's state equation
will turn into a standard backward stochastic Riccati equation (see Tang
\cite{Tang03}). Moreover, under certain assumptions, the stochastic Riccati
equation admits a unique solution. For a general case, namely, \emph{%
non-convex} control set, see Section \ref{sec4}.

\item We focus on the LQ Stackelberg game with the control constrained in a
closed convex set $\Gamma $ of full space: $\Gamma \subset \mathbb{R}^{m}$.
One of the motivations to study the LQ problems with control constraint
arises naturally from mathematical finance. For instance, the no-shorting
constraint\footnote{%
Short sales have so many risks that make it unsuitable for a novice
investor. For starters, ~there is theoretically no limit to the investor's
possible loss if the stock price rises instead of declines. A stock can only
fall to zero, resulting in a $100\%$ loss for a long investor, but there is
no limit to how high a stock can theoretically go. A short seller who has
not covered his or her position with a stop-loss buyback order can suffer
tremendous losses if the stock price runs higher. For example, consider a
company that becomes embroiled in scandal when its stock is trading at $\$70$
per share. An investor sees an opportunity to make a quick profit and sells
the stock short at $\$65$. But then the company is able to quickly exonerate
itself from the accusations by coming up with tangible proof to the
contrary. The stock price quickly rises to $\$80$ a share, leaving the
investor with a loss of $\$15$ per share for the moment. If the stock
continues to rise, so do the investor's losses. Besides, short selling also
involves significant expenses. There are the costs of borrowing the security
to sell, the interest payable on the margin account that holds it, and
trading commissions, etc.} in portfolio selection leads to the LQ control
with positive control ($\Gamma =\mathbb{R}_{+}^{m},$ the positive orthant).
Moreover, since the general market accessibility constraint, it also
promises interesting to investigate the LQ control with more general closed
convex cone constraint (see \cite{hz}). As a response, this paper
investigates the LQ Stackelberg game with general closed convex control
constraint. The control constraint will bring some new features here: (1)
The related Hamitonian system is no longer linear, and it becomes a class of
nonlinear FBSDEs with projection operator. (2) Due to the nonlinearity, the
standard Riccati equation with feedback control is no longer valid to
represent the open-loop solution to the two-person leader-follower
stochastic differential game.
\end{itemize}

The rest of the paper is organized as follows. Section \ref{sec2} is devoted
to presenting the maximum principle for a Stackelberg game of follower and
leader under the AOL information pattern with convex control input, which is
well known (cf. \cite{Wu1998}). Based on previous result, we study
Stackelberg games under the ACLM information pattern, and establish the
maximum principle for the leader's optimal strategy, together with some
other preliminary results. In section \ref{sec3}, as applications, linear
quadratic Stackelberg games under the AOL and ACLM information patterns are
investigated, respectively. For former case, on the one hand, we prove the
the existence and uniqueness of the solution to the associated Hamiltonian
system for follower with projection operator; On the other hand, we show the
existence and uniqueness of the solution to the associated stochastic
Riccati equation under certain assumptions. For the latter case, we merely
derive the associated Riccati equation which consists of a kind of complex
FBSDEs, due to the quadratic and irregular feature. Some conclusions and
unsolved issues for future research are displayed in Section \ref{sec4}.
Some proof and discussion are displayed in Appendix.

\section{Preliminaries}

\label{sec2}

Throughout this paper, we denote by $\mathbb{R}^{n}$ the space of $n$%
-dimensional Euclidean space, by $\mathbb{R}^{n\times d}$ the space the
matrices with order $n\times d$. Let $(\Omega ,\mathcal{F},\{\mathcal{F}%
_{t}\}_{t\geq 0},P)$ be a complete filtered probability space on which a $1$%
-dimensional standard Brownian motion $W(\cdot )$ is defined, with $\{%
\mathcal{F}_{t}\}_{t\geq 0}$ being its natural filtration, augmented by all
the $P$-null sets.

We now introduce the following spaces of processes:
\begin{align*}
\mathcal{S}^{2}(0,T;\mathbb{R}^{n})& \triangleq \left\{ \mathbb{R}^{n}\text{%
-valued }\mathcal{F}_{t}\text{-adapted process }\phi (t)\text{; }\mathbb{E}%
\left[ \sup\limits_{0\leq t\leq T}\left\vert \phi _{t}\right\vert ^{2}\right]
<\infty \right\} , \\
\mathcal{M}^{2}(0,T;\mathbb{R}^{n})& \triangleq \left\{ \mathbb{R}^{n}\text{%
-valued }\mathcal{F}_{t}\text{-adapted process }\varphi (t)\text{; }\mathbb{E%
}\left[ \int_{0}^{T}\left\vert \varphi _{t}\right\vert ^{2}\mbox{\rm d}t%
\right] <\infty \right\} ,
\end{align*}%
and denote $\mathcal{N}^{2}\left[ 0,T\right] =\mathcal{S}^{2}(0,T;\mathbb{R}%
^{n})\times \mathcal{S}^{2}(0,T;\mathbb{R}^{n})\times \mathcal{M}^{2}(0,T;%
\mathbb{R}^{n}).$ Clearly, $\mathcal{N}^{2}\left[ 0,T\right] $ forms a
Banach space.

Consider the following:
\begin{equation}
\left\{
\begin{array}{rcl}
\mathrm{d}x\left( t\right) & = & b\left( t,x\left( t\right) ,u\left(
t\right) ,v\left( t\right) \right) \mathrm{d}t+\sigma \left( t,x\left(
t\right) ,u\left( t\right) ,v\left( t\right) \right) \mathrm{d}W\left(
t\right) , \\
x\left( 0\right) & = & x_{0}\in \mathbb{R}^{n},%
\end{array}%
\right.  \label{sde1}
\end{equation}%
where and $\left( u\left( \cdot \right) ,v\left( \cdot \right) \right) $
denotes the decisions of the \textit{leader} and the \textit{follower}, with
values in subsets $U$ and $V$ in some closed convex subset $\Gamma _{1}$ and
$\Gamma _{2}$ of full space $\mathbb{R}^{m_{1}}$ and $\mathbb{R}^{m_{2}},$
respectively.

The cost functionals for the leader and the follower to minimize are given,
respectively, as follows:%
\begin{equation}
\mathcal{J}_{1}\left( u,v\right) =\mathbb{E}\left[ \int_{0}^{T}l_{1}\left(
t,x\left( t\right) ,u\left( t\right) ,v\left( t\right) \right) \mathrm{d}%
t+\Phi _{1}\left( x\left( T\right) \right) \right]  \label{c1}
\end{equation}%
and
\begin{equation}
\mathcal{J}_{2}\left( u,v\right) =\mathbb{E}\left[ \int_{0}^{T}l_{2}\left(
t,x\left( t\right) ,u\left( t\right) ,v\left( t\right) \right) \mathrm{d}%
t+\Phi _{2}\left( x\left( T\right) \right) \right] .  \label{c2}
\end{equation}%
The coefficients $b$ and $\sigma $ in (\ref{sde1}), and $l_{i}$ and $\Phi
_{i}$, $i=1,2$ in the cost functionals (\ref{c1}) and (\ref{c2}) are
specified as follows:%
\begin{eqnarray*}
b &:&\Omega \times \left[ 0,T\right] \times \mathbb{R}^{n}\times \mathbb{R}%
^{m_{1}}\times \mathbb{R}^{m_{2}}\rightarrow \mathbb{R}^{n},\text{ }\mathcal{%
P\times B}\left( \mathbb{R}^{n+m_{1}+m_{2}}\right) /\mathcal{B}\left(
\mathbb{R}^{n}\right) \text{ measurable,} \\
\sigma &:&\Omega \times \left[ 0,T\right] \times \mathbb{R}^{n}\times
\mathbb{R}^{m_{1}}\times \mathbb{R}^{m_{2}}\rightarrow \mathbb{R}^{n\times
d},\text{ }\mathcal{P\times B}\left( \mathbb{R}^{n+m_{1}+m_{2}}\right) /%
\mathcal{B}\left( \mathbb{R}^{n}\right) \text{ measurable,} \\
l_{i} &:&\Omega \times \left[ 0,T\right] \times \mathbb{R}^{n}\times \mathbb{%
R}^{m_{1}}\times \mathbb{R}^{m_{2}}\rightarrow \mathbb{R},\text{ }\mathcal{%
P\times B}\left( \mathbb{R}^{n}\right) \times \mathcal{B}\left( \Gamma
_{1}\right) \times \mathcal{B}\left( \Gamma _{2}\right) /\mathcal{B}\left(
\mathbb{R}\right) \text{ measurable,} \\
\Phi _{i} &:&\Omega \times \mathbb{R}^{n}\rightarrow \mathbb{R},\text{ }%
\mathcal{F}_{T}\mathcal{\times B}\left( \mathbb{R}^{n}\right) /\mathcal{B}%
\left( \mathbb{R}\right) \text{ measurable.}
\end{eqnarray*}%
Letting $\varphi \left( t,x,u,v\right) =b\left( t,x,u,v\right) ,$ $\sigma
\left( t,x,u,v\right) ,$ $l_{i}\left( t,x,u,v\right) ,$ $\Phi _{i},$ $i=1,2,$
we give the standing assumptions of our paper:

\begin{description}
\item[(A1)] We postulate throughout the paper that $\varphi $ and its first
and second derivatives are uniformly Lipschitz with respect to $\left(
x,u,v\right) $ and $\varphi \left( \cdot ,x,u,v\right) \in \mathcal{M}^{2}$,
for $\left( x,u,v\right) \in \mathbb{R}^{n}\times \mathbb{R}^{m_{1}}\times
\mathbb{R}^{m_{2}}.$
\end{description}

In this paper, we focus on the players' information structures $\eta $ as
follows:

(a) adapted open-loop (AOL): $\eta \left( t\right) =\left\{ x_{0},\mathcal{F}%
_{t}\right\} $, $t\in \left[ 0,T\right] ,$

(b) adapted closed-loop memoryless (ACLM): $\eta \left( t\right) =\left\{
x_{0},x\left( t\right) ,\mathcal{F}_{t}\right\} $, $t\in \left[ 0,T\right] .$

\subsection{AOL information structure}

For the AOL information structure, the admissible strategy spaces for the
leader and the follower are denoted by
\begin{equation*}
\mathcal{U}=\left\{ u\left\vert u:\Omega \times \left[ 0,T\right]
\rightarrow \Gamma _{1}\text{ is }\mathcal{F}_{t}\text{-adapted satisfying }%
\mathbb{E}\left[ \int_{0}^{T}\left\vert u\left( t\right) \right\vert ^{2}%
\mathrm{d}t\right] <\infty \right. \right\} ,
\end{equation*}%
\begin{equation*}
\mathcal{V}=\bigg \{v\bigg |v:\Omega \times \left[ 0,T\right] \times
\mathcal{U}\rightarrow \Gamma _{2},\text{ for }u\in \mathcal{U},\text{ }%
v\left( \cdot ,u\right) \text{ is }\mathcal{F}_{t}\text{-adapted satisfying }%
\mathbb{E}\left[ \int_{0}^{T}\left\vert v\left( t\right) \right\vert ^{2}%
\mathrm{d}t\right] <\infty \bigg \}.
\end{equation*}

\begin{remark}
Since the initial state $x_{0}$ is commonly known by both players, $x_{0}$
is suppressed.
\end{remark}

Under the AOL information pattern, we first establish a maximum principle
for optimal control of the follower, whenever given the leader's strategy $%
u\in \mathcal{U}$.

\noindent \textbf{Problem (AOL-F)} Fix $u\in \mathcal{U}.$ Seek an
admissible control $v^{\ast }\left( \cdot \right) \in \mathcal{V}$ such that
\begin{equation*}
\mathcal{J}_{2}\left( u,v^{\ast }\right) =\inf_{v\left( \cdot \right) \in
\mathcal{V}}\mathcal{J}_{2}\left( u,v\right)
\end{equation*}%
subject to (\ref{sde1}).

Define the Hamiltonian function:%
\begin{eqnarray*}
\mathcal{H}_{2}\left( t,x,u,v,p_{2},q_{2}\right) &=&\left\langle
p_{2},b\left( t,x,u,v\right) \right\rangle +\left\langle q_{2},\sigma \left(
t,x,u,v\right) \right\rangle +l_{2}\left( t,x,u,v\right) , \\
\forall \left( t,x,u,v,p_{2},q_{2}\right) &\in &\left[ 0,T\right] \times
\mathbb{R}^{n}\times \mathbb{R}^{m_{1}}\times \mathbb{R}^{m_{2}}\times
\mathbb{R}^{n}\times \mathbb{R}^{n}.
\end{eqnarray*}%
Then the maximum principle\footnote{%
Since the control region is closed and convex, the first-order adjoint
equation is needed. For general case, that is, compact control domain, the
second-order adjoint equation must be introduced (see \cite{YZ} and \cite%
{Peng1990} for more details).} (cf. \cite{YZ}) says that if we assume that $%
v^{\ast }\left( \cdot \right) \in \mathcal{V}$ is an optimal control, there
exists a unique adapted solution $\left( p_{2}\left( \cdot \right)
,q_{2}\left( \cdot \right) \right) \in \mathcal{S}^{2}(0,T;\mathbb{R}%
^{n})\times \mathcal{M}^{2}(0,T;\mathbb{R}^{n\times d})$ such that%
\begin{equation}
\left\{
\begin{array}{rcl}
\mathrm{d}x\left( t\right) & = & b\left( t,x\left( t\right) ,u\left(
t\right) ,v^{\ast }\left( t\right) \right) \mathrm{d}t+\sigma \left(
t,x\left( t\right) ,u\left( t\right) ,v^{\ast }\left( t\right) \right)
\mathrm{d}W\left( t\right) , \\
-\mathrm{d}p_{2}\left( t\right) & = & \frac{\partial }{\partial x}\mathcal{H}%
_{2}\left( t,x\left( t\right) ,u\left( t\right) ,v^{\ast }\left( t\right)
,p_{2}\left( t\right) ,q_{2}\left( t\right) \right) \mathrm{d}t-q_{2}\left(
t\right) \mathrm{d}W\left( t\right) , \\
x\left( 0\right) & = & x_{0},\text{ }p_{2}\left( T\right) =\frac{\partial }{%
\partial x}\Phi _{2}\left( x\left( T\right) \right) ,%
\end{array}%
\right.  \label{a1}
\end{equation}%
and
\begin{equation}
v^{\ast }\left( t,x\left( t\right) ,u\left( t\right) ,p_{2}\left( t\right)
,q_{2}\left( t\right) \right) =\arg \min_{v\in \Gamma _{2}}\mathcal{H}%
_{2}\left( t,x\left( t\right) ,u\left( t\right) ,v,p_{2}\left( t\right)
,q_{2}\left( t\right) \right) .  \label{o1}
\end{equation}%
We suppose that, the function $v^{\ast }\left( t,x,u,p_{2},q_{2}\right) $ is
uniquely defined and is uniformly Lipschitz continuous with respect to $%
\left( x,u,p_{2},q_{2}\right) $ and continuously differentiable\footnote{%
In fact, we will see later in Section \ref{sec3}, in the framework of linear
quadratic, whenever the control domain is closed and convex, by convex
analysis, there indeed exists a unique optimal control in the form of
projector.}. Now inserting $v^{\ast }\left( t,x,u,p_{2},q_{2}\right) $ into (%
\ref{a1}), we formulate the optimal control problem for leader.

\noindent \textbf{Problem (AOL-L)} Seek an admissible control $u^{\ast
}\left( \cdot \right) \in \mathcal{U}$ such that
\begin{equation*}
\mathcal{J}_{1}\left( u^{\ast },v^{\ast }\right) =\inf_{u\left( \cdot
\right) \in \mathcal{U}}\mathcal{J}_{1}\left( u,v^{\ast }\right)
\end{equation*}%
subject to
\begin{equation}
\left\{
\begin{array}{rcl}
\mathrm{d}x\left( t\right) & = & b\left( t,x\left( t\right) ,u\left(
t\right) ,v^{\ast }\left( t,x\left( t\right) ,u\left( t\right) ,p_{2}\left(
t\right) ,q_{2}\left( t\right) \right) \right) \mathrm{d}t \\
&  & +\sigma \left( t,x\left( t\right) ,u\left( t\right) ,v^{\ast }\left(
t,x\left( t\right) ,u\left( t\right) ,p_{2}\left( t\right) ,q_{2}\left(
t\right) \right) \right) \mathrm{d}W\left( t\right) , \\
-\mathrm{d}p_{2}\left( t\right) & = & \frac{\partial }{\partial x}\mathcal{H}%
_{2}\left( t,x\left( t\right) ,u\left( t\right) ,v^{\ast }\left( t,x\left(
t\right) ,u\left( t\right) ,p_{2}\left( t\right) ,q_{2}\left( t\right)
\right) ,p_{2}\left( t\right) ,q_{2}\left( t\right) \right) \mathrm{d}t \\
&  & -q_{2}\left( t\right) \mathrm{d}W\left( t\right) , \\
x\left( 0\right) & = & x_{0},\text{ }p_{2}\left( T\right) =\frac{\partial }{%
\partial x}\Phi _{2}\left( x\left( T\right) \right) .%
\end{array}%
\right.  \label{a2}
\end{equation}

\begin{remark}
From Problem (AOL-L), reader maybe realize that one of the motivations to
focus on fully coupled FBSDEs comes from stochastic Stackelberg differential
games. In the literature, to treat fully coupled FBSDEs, Antonelli first
provided a counterexample (see \cite{Anto1993}) showing that the Lipschitz
condition is not enough for the existence of FBSDEs in an arbitrarily large
time duration. Thereby, more assumptions are essentially needed. To the best
of our knowledge, there exist two approaches to deal with such FBSDEs. The
first one is purely probabilistic (see \cite{Anto1993,hp, PW1999}) under the
monotone conditions; The second one concerns a kind of so called
\textquotedblleft four-steps scheme\textquotedblright\ method (combination
of the methods of partial differential equation and probability or
stochastic optimal control). Several major applications in mathematical
finance have been made (see \cite{DE, PW1999}). It is now very clear that
certain important problems in mathematical economics and mathematical
finance, especially in the optimization problem, are formulated to be fully
coupled FBSDEs (see \cite{WY2014}).
\end{remark}

Clearly, for given $u\left( \cdot \right) \in \mathcal{U},$ FBSDEs (\ref{a2}%
) are fully coupled (while FBSDEs (3.2) in Bensoussan et al. \cite{BCS2015}
are partially coupled). In order to make the leader's problem well-posed, we
proposed some conditions, mainly taken from Hu and Peng \cite{hp} (see also
Peng and Wu \cite{PW1999}) to guarantee that FBSDEs (\ref{a2}) admit a
unique adapted solution.

For $\Lambda ^{1}=\left( x^{1},p_{2}^{1},q_{2}^{1}\right) \in \mathbb{R}%
^{n}\times \mathbb{R}^{n}\times \mathbb{R}^{n},$ $\Lambda ^{2}=\left(
x^{2},p_{2}^{2},q_{2}^{2}\right) \in \mathbb{R}^{n}\times \mathbb{R}%
^{n}\times \mathbb{R}^{n},$%
\begin{equation*}
\left[ \Lambda ^{1},\Lambda ^{2}\right] =\left\langle
x^{1},x^{2}\right\rangle +\left\langle p_{2}^{1},p_{2}^{2}\right\rangle
+\left\langle q_{2}^{1},q_{2}^{2}\right\rangle ,
\end{equation*}%
Let
\begin{equation*}
f_{2}\left( t,u,x,p_{2},q_{2}\right) =\frac{\partial }{\partial x}\mathcal{H}%
_{2}\left( t,x,u,v^{\ast }\left( t,x,u,p_{2},q_{2}\right) ,p_{2},q_{2}\right)
\end{equation*}%
For $\Lambda =\left( x,p_{2},q_{2}\right) \in \mathbb{R}^{n}\times \mathbb{R}%
^{n}\times \mathbb{R}^{n},$%
\begin{equation*}
F\left( t,u,\Lambda \right) =\left( f_{2}\left( t,u,\Lambda \right) ,b\left(
t,u,\Lambda \right) ,\sigma \left( t,u,\Lambda \right) \right) .
\end{equation*}%
We assume

\begin{description}
\item[(A2)] For each $\Lambda =\left( x,p_{2},q_{2}\right) \in \mathbb{R}%
^{n}\times \mathbb{R}^{n}\times \mathbb{R}^{n},$ $F\left( t,u,\Lambda
\right) \in \mathcal{M}^{2}\left( 0,T;\mathbb{R}^{n}\times \mathbb{R}%
^{n}\times \mathbb{R}^{n}\right) ,$ and for each $x\in \mathbb{R}^{n},$ $%
\frac{\partial }{\partial x}\Phi _{2}\left( x\right) \in L^{2}\left( \Omega ,%
\mathcal{F}_{T};\mathbb{R}^{n}\right) ;$ there exists a constant $c_{1}>0,$
such that
\begin{equation*}
\left\vert F\left( t,u,\Lambda ^{1}\right) -F\left( t,u,\Lambda ^{2}\right)
\right\vert \leq c_{1}\left\vert \Lambda ^{1}-\Lambda ^{2}\right\vert ,\text{
}\Lambda ^{i}\in \mathbb{R}^{n}\times \mathbb{R}^{n}\times \mathbb{R}^{n},%
\text{ }i=1,2
\end{equation*}%
and
\begin{equation*}
\left\vert \frac{\partial }{\partial x}\Phi _{2}\left( x_{1}\right) -\frac{%
\partial }{\partial x}\Phi _{2}\left( x_{2}\right) \right\vert \leq
c_{1}\left\vert x_{1}-x_{2}\right\vert ,\text{ }P\text{-a.s. }x_{1},x_{2}\in
\mathbb{R}^{n}.
\end{equation*}

\item[(A3)] There exists a constant $c_{2}>0,$ such that
\begin{equation*}
\left[ F\left( t,u,\Lambda ^{1}\right) -F\left( t,u,\Lambda ^{2}\right)
,\Lambda ^{1}-\Lambda ^{2}\right] \leq -c_{2}\left\vert \Lambda ^{1}-\Lambda
^{2}\right\vert ^{2},\Lambda ^{i}\in \mathbb{R}^{n}\times \mathbb{R}%
^{n}\times \mathbb{R}^{n},\text{ }i=1,2
\end{equation*}%
and
\begin{equation*}
\left\langle \frac{\partial }{\partial x}\Phi _{2}\left( x_{1}\right) -\frac{%
\partial }{\partial x}\Phi _{2}\left( x_{2}\right) ,x_{1}-x_{2}\right\rangle
\geq c_{2}\left\vert x_{1}-x_{2}\right\vert ,\text{ }P\text{-a.s. }%
x_{1},x_{2}\in \mathbb{R}^{n}.
\end{equation*}
\end{description}

Under (A2)-(A3), FBSDEs (\ref{a2}) have a unique adapted solution for $u\in
\mathcal{U}$. Next, we will establish a maximum principle for leader's
optimal control. Since the control domain is convex, the first adjoint
equation is needed.

\begin{proposition}
\label{pro1}Suppose that the Assumptions \emph{(A2)-(A3)} hold. Let $u^{\ast
}\left( \cdot \right) $ is an optimal strategy for the leader. Then there
exists a unique adapted solution\newline
$\left( k\left( \cdot \right) ,p_{1}\left( \cdot \right) ,q_{1}\left( \cdot
\right) \right) \in \mathcal{M}^{2}\left( 0,T;\mathbb{R}^{n}\times \mathbb{R}%
^{n}\times \mathbb{R}^{n}\right) ,$ such that%
\begin{equation*}
u^{\ast }\left( t\right) =\arg \min_{u\in \mathcal{U}}\mathcal{H}_{1}\left(
t,u,x\left( t\right) ,k\left( t\right) ,p_{1}\left( \cdot \right)
,p_{2}\left( \cdot \right) ,q_{1}\left( \cdot \right) ,q_{2}\left( \cdot
\right) \right) ,
\end{equation*}%
where
\begin{eqnarray*}
\mathcal{H}_{1}\left( t,u,x,k,p_{1},p_{2},q_{1},q_{2}\right) &=&\left\langle
p_{1},b\left( t,x,u,v^{\ast }\left( t,x,u,p_{2},q_{2}\right) \right)
\right\rangle \\
&&\left\langle q_{1},\sigma \left( t,x,u,v^{\ast }\left(
t,x,u,p_{2},q_{2}\right) \right) \right\rangle \\
&&+l_{1}\left( t,x,u,v^{\ast }\left( t,x,u,p_{2},q_{2}\right) \right) \\
&&-\left\langle k,f_{2}\left( t,u,x,p_{2},q_{2}\right) \right\rangle
\end{eqnarray*}%
and
\begin{equation}
\left\{
\begin{array}{rcl}
\mathrm{d}k\left( t\right) & = & -\frac{\partial }{\partial p_{2}}\mathcal{H}%
_{1}\mathrm{d}t-\frac{\partial }{\partial q_{2}}\mathcal{H}_{1}\mathrm{d}%
W\left( t\right), \\
\mathrm{d}p_{1}\left( t\right) & = & \frac{\partial }{\partial x}\mathcal{H}%
_{1}\mathrm{d}t+q_{1}\left( t\right) \mathrm{d}W\left( t\right) , \\
k\left( 0\right) & = & 0,\text{ }p_{1}\left( T\right) =-\frac{\partial ^{2}}{%
\partial x^{2}}\Phi _{2}\left( x\left( T\right) \right) k\left( T\right) +%
\frac{\partial }{\partial x}\Phi _{1}\left( x\left( T\right) \right) ,%
\end{array}%
\right.  \label{a3}
\end{equation}%
where
\begin{eqnarray*}
\frac{\partial }{\partial p_{2}}\mathcal{H}_{1} &=&\left( \frac{\partial b}{%
\partial v}\frac{\partial v^{\ast }}{\partial p_{2}}\right) ^{\top
}p_{1}+\left( \frac{\partial \sigma }{\partial v}\frac{\partial v^{\ast }}{%
\partial p_{2}}\right) ^{\top }q_{1}+\left( \frac{\partial v^{\ast }}{%
\partial p_{2}}\right) ^{\top }\frac{\partial l_{1}}{\partial v} \\
&&-\frac{\partial b}{\partial x}k-\sum_{i=1}^{n}k_{i}\left( \frac{\partial
v^{\ast }}{\partial p_{2}}\right) ^{\top }\frac{\partial }{\partial v}\left(
\frac{\partial b}{\partial x_{i}}\right) ^{\top }p_{2} \\
&&-\left( \frac{\partial ^{2}l_{2}}{\partial v\partial x}\frac{\partial
v^{\ast }}{\partial p_{2}}\right) ^{\top }k-\sum_{i=1}^{n}k_{i}\left( \frac{%
\partial v^{\ast }}{\partial p_{2}}\right) ^{\top }\frac{\partial }{\partial
v}\left( \frac{\partial \sigma }{\partial x_{i}}\right) ^{\top }q_{2}, \\
\frac{\partial }{\partial q_{2}}\mathcal{H}_{1} &=&\left( \frac{\partial b}{%
\partial v}\frac{\partial v^{\ast }}{\partial q_{2}}\right) ^{\top
}p_{1}+\left( \frac{\partial \sigma }{\partial v}\frac{\partial v^{\ast }}{%
\partial q_{2}}\right) ^{\top }q_{1}+\left( \frac{\partial v^{\ast }}{%
\partial q_{2}}\right) ^{\top }\frac{\partial l_{1}}{\partial v} \\
&&-\sum_{i=1}^{n}k_{i}\left( \frac{\partial v^{\ast }}{\partial q_{2}}%
\right) ^{\top }\frac{\partial }{\partial v}\left( \frac{\partial b}{%
\partial x_{i}}\right) ^{\top }p_{2} \\
&&-\frac{\partial \sigma }{\partial x}k-\sum_{i=1}^{n}k_{i}\left( \frac{%
\partial v^{\ast }}{\partial q_{2}}\right) ^{\top }\frac{\partial }{\partial
v}\left( \frac{\partial \sigma }{\partial x_{i}}\right) ^{\top }q_{2}-\left(
\frac{\partial ^{2}l_{2}}{\partial v\partial x}\frac{\partial v^{\ast }}{%
\partial q_{2}}\right) ^{\top }k, \\
\frac{\partial }{\partial x}\mathcal{H}_{1} &=&\left( \frac{\partial b}{%
\partial x}\right) ^{\top }p_{1}+\left( \frac{\partial b}{\partial v}\frac{%
\partial v^{\ast }}{\partial x}\right) ^{\top }p_{1}+\left( \frac{\partial
\sigma }{\partial x}\right) ^{\top }q_{1}+\left( \frac{\partial \sigma }{%
\partial v}\frac{\partial v^{\ast }}{\partial x}\right) ^{\top }q_{1} \\
&&+\frac{\partial l_{1}}{\partial x}+\left( \frac{\partial v^{\ast }}{%
\partial x}\right) ^{\top }\frac{\partial l_{1}}{\partial v} \\
&&-\sum_{i=1}^{n}k_{i}\left[ \frac{\partial }{\partial x}\left( \frac{%
\partial b}{\partial x_{i}}\right) ^{\top }+\left( \frac{\partial v^{\ast }}{%
\partial x_{2}}\right) ^{\top }\frac{\partial }{\partial v}\left( \frac{%
\partial b}{\partial x_{i}}\right) ^{\top }\right] p_{2} \\
&&-\sum_{i=1}^{n}k_{i}\left[ \frac{\partial }{\partial x}\left( \frac{%
\partial \sigma }{\partial x_{i}}\right) ^{\top }+\left( \frac{\partial
v^{\ast }}{\partial x_{2}}\right) ^{\top }\frac{\partial }{\partial v}\left(
\frac{\partial \sigma }{\partial x_{i}}\right) ^{\top }\right] q_{2} \\
&&-\left( \frac{\partial ^{2}l_{2}}{\partial x^{2}}+\frac{\partial ^{2}l_{2}%
}{\partial x\partial v}\frac{\partial v^{\ast }}{\partial x}\right) ^{\top
}k.
\end{eqnarray*}
\end{proposition}


\subsection{ACLM information structure}

For the ACLM information structure, the admissible strategy spaces for the
leader and the follower are denoted by
\begin{eqnarray*}
\mathcal{U} &=&\Bigg \{u|u:\Omega \times \left[ 0,T\right] \times \mathbb{R}%
^{n}\rightarrow \Gamma _{1}\text{ is }\mathcal{F}_{t}\text{-adapted for any }%
x\in \mathbb{R}^{n},\text{ }u\left( t,x\right) \text{ is continuously} \\
&&\text{differentiable in }x\text{ for any }\left( t,x\right) \in \Omega
\times \left[ 0,T\right] \text{ satisfying }\left\vert \frac{\partial u}{%
\partial x}\right\vert <K\Bigg \}, \\
\mathcal{V} &=&\left\{ v\left\vert :\Omega \times \left[ 0,T\right] \times
\mathbb{R}^{n}\times \mathcal{U}\rightarrow \Gamma _{2}\text{ is }\mathcal{F}%
_{t}\text{-adapted for any }x\in \mathbb{R}^{n}\text{ and }u\in \mathcal{U}%
\right. \right\} .
\end{eqnarray*}%
Now consider the following,

\noindent \textbf{Problem (ACLM-F)} For any $u\in \mathcal{U}$, seek an
admissible control $v^{\ast }\left( \cdot \right) \in \mathcal{V}$ such that
\begin{equation*}
\mathcal{J}_{2}\left( u,v^{\ast }\right) =\inf_{v\left( \cdot \right) \in
\mathcal{V}}\mathcal{J}_{2}\left( u,v\right)
\end{equation*}%
subject to
\begin{equation}
\left\{
\begin{array}{rcl}
\mathrm{d}x\left( t\right) & = & b\left( t,x\left( t\right) ,u\left(
t,x\left( t\right) \right) ,v\right) \mathrm{d}t+\sigma \left( t,x\left(
t\right) ,u\left( t,x\left( t\right) \right) ,v\right) \mathrm{d}W\left(
t\right) , \\
x\left( 0\right) & = & x_{0}.%
\end{array}%
\right.  \label{sde3}
\end{equation}%
Once again, then the maximum principle (cf. \cite{YZ}) states that if we
assume that $v^{\ast }\left( \cdot \right) \in \mathcal{V}$ is an optimal
control, there exists a unique adapted solution $\left( p_{2}\left( \cdot
\right) ,q_{2}\left( \cdot \right) \right) \in \mathcal{S}^{2}(0,T;\mathbb{R}%
^{n})\times \mathcal{M}^{2}(0,T;\mathbb{R}^{n\times d})$ such that%
\begin{equation}
\left\{
\begin{array}{rcl}
-\mathrm{d}p_{2}\left( t\right) & = & \frac{\partial }{\partial x}\mathcal{H}%
_{2}\left( t,x\left( t\right) ,u\left( t,x\left( t\right) \right) ,v\left(
t\right) ,p_{2}\left( t\right) ,q_{2}\left( t\right) \right) \mathrm{d}%
t-q_{2}\left( t\right) \mathrm{d}W\left( t\right) , \\
p_{2}\left( T\right) & = & \frac{\partial }{\partial x}\Phi _{2}\left(
x\left( T\right) \right) ,%
\end{array}%
\right.  \label{a4}
\end{equation}%
with
\begin{eqnarray*}
&&\frac{\partial }{\partial x}\mathcal{H}_{2}\left( t,x\left( t\right)
,u\left( t,x\left( t\right) \right) ,v\left( t\right) ,p_{2}\left( t\right)
,q_{2}\left( t\right) \right) \\
&=&\left( \frac{\partial b}{\partial x}+\frac{\partial b}{\partial u}\frac{%
\partial u}{\partial x}\right) ^{\top }p_{2}+\left( \frac{\partial \sigma }{%
\partial x}+\frac{\partial \sigma }{\partial u}\frac{\partial u}{\partial x}%
\right) ^{\top }q_{2} \\
&&+\frac{\partial l_{2}}{\partial x}+\left( \frac{\partial u}{\partial x}%
\right) ^{\top }\frac{\partial l_{2}}{\partial u}
\end{eqnarray*}%
and
\begin{equation}
v^{\ast }\left( t,x\left( t\right) ,u,p_{2}\left( t\right) ,q_{2}\left(
t\right) \right) =\arg \min_{v\in \Gamma _{2}}\mathcal{H}_{2}\left(
t,x\left( t\right) ,u\left( t,x\left( t\right) \right) ,v,p_{2}\left(
t\right) ,q_{2}\left( t\right) \right) .  \label{ACLM-f}
\end{equation}%
We now postulate that for any leader's strategy $u\in \mathcal{U}$, there
exists a unique strategy $v^{\ast }$ for the follower that minimizes his
cost functional $\mathcal{J}_{2}$ and that (\ref{ACLM-f}) yields $v^{\ast
}\left( t\right) =v^{\ast }\left( t,x\left( t\right) ,u,p_{2}\left( t\right)
,q_{2}\left( t\right) \right) $. Then, taking into account the follower's
optimal response, the leader will be taken action by solving the optimal
control problem:
\begin{equation}
\mathcal{J}_{1}\left( u^{\ast },v^{\ast }\right) =\inf_{v\left( \cdot
\right) \in \mathcal{V}}\mathcal{J}_{1}\left( u,v^{\ast }\right)
\label{cost2}
\end{equation}%
subject to
\begin{equation}
\left\{
\begin{array}{rcl}
\mathrm{d}x\left( t\right) & = & b\left( t,x\left( t\right) ,u\left(
t,x\left( t\right) \right) ,v^{\ast }\left( t\right) \right) \mathrm{d}%
t+\sigma \left( t,x\left( t\right) ,u\left( t,x\left( t\right) \right)
,v^{\ast }\left( t\right) \right) \mathrm{d}W\left( t\right) , \\
-\mathrm{d}p_{2}\left( t\right) & = & \frac{\partial }{\partial x}\mathcal{H}%
_{2}\left( t,x\left( t\right) ,u\left( t\right) ,v^{\ast }\left( t\right)
,p_{2}\left( t\right) ,q_{2}\left( t\right) \right) \mathrm{d}t-q_{2}\left(
t\right) \mathrm{d}W\left( t\right) , \\
x\left( 0\right) & = & x_{0},\text{ }p_{2}\left( T\right) =\frac{\partial }{%
\partial x}\Phi _{2}\left( x\left( T\right) \right) ,%
\end{array}%
\right.  \label{a5}
\end{equation}%
with
\begin{eqnarray*}
&&\frac{\partial }{\partial x}\mathcal{H}_{2}\left( t,x\left( t\right)
,u\left( t,x\left( t\right) \right) ,v^{\ast }\left( t\right) ,p_{2}\left(
t\right) ,q_{2}\left( t\right) \right) \\
&=&\left( \frac{\partial b}{\partial x}+\frac{\partial b}{\partial u}\frac{%
\partial u}{\partial x}\right) ^{\top }p_{2}+\left( \frac{\partial \sigma }{%
\partial x}+\frac{\partial \sigma }{\partial u}\frac{\partial u}{\partial x}%
\right) ^{\top }q_{2} \\
&&+\frac{\partial l_{2}}{\partial x}+\left( \frac{\partial u}{\partial x}%
\right) ^{\top }\frac{\partial l_{2}}{\partial u}.
\end{eqnarray*}

We assume that the follower has a unique optimal response strategy $v^{\ast
} $ for every strategy $u\in \mathcal{U}$ of the leader. Similarly, we
further suppose that the leader's problem is well-posed, i.e., for each $%
u\in \mathcal{U}$, there exists a unique triple $\left( x\left( \cdot
\right) ,p_{2}\left( \cdot \right) ,q_{2}\left( \cdot \right) \right) \in
\mathcal{N}^{2}\left[ 0,T\right] $ that solves FBSDEs (\ref{a5}). Clearly,
the appearance of the derivative $\frac{\partial u}{\partial x}$ of the
strategy $u$ in (\ref{a5}) leads to in a non-standard optimal control
problem for the leader. Employing the idea from Bensoussan et al. \cite%
{BCS2015}, we first transform the original issue to a standard stochastic
optimal control problem, and subsequently establish the equivalence between
the two in the sense that they coincide with the same optimal trajectory and
cost.

Let us introduce the following optimal control problem:%
\begin{eqnarray}
\mathcal{J}\left( u_{1}^{\ast },u_{2}^{\ast }\right) &=&\min_{\left(
u_{1}\left( \cdot \right) ,u_{2}\left( \cdot \right) \right) }\mathcal{J}%
_{1}\left( u,v\right)  \notag \\
&=&\mathbb{E}\Bigg [\int_{0}^{T}l_{1}\left( t,x\left( t\right) ,u_{2}\left(
t\right) x\left( t\right) +u_{1}\left( t\right) ,\mu ^{\ast }\left( t\right)
\right) \mathrm{d}t+\Phi _{1}\left( x\left( T\right) \right) \Bigg ],
\label{necost}
\end{eqnarray}%
where
\begin{equation*}
\mu ^{\ast }\left( t\right) =v^{\ast }\left( t,x\left( t\right) ,u_{2}\left(
t\right) x\left( t\right) +u_{1}\left( t\right) ,p_{2}\left( t\right)
,q_{2}\left( t\right) \right) ,
\end{equation*}%
subject to
\begin{equation}
\left\{
\begin{array}{rcl}
\mathrm{d}x\left( t\right) & = & b\left( t,x\left( t\right) ,u_{2}\left(
t\right) x\left( t\right) +u_{1}\left( t\right) ,\mu ^{\ast }\left( t\right)
\right) \mathrm{d}t \\
&  & +\sigma \left( t,x\left( t\right) ,u_{2}\left( t\right) x\left(
t\right) +u_{1}\left( t\right) ,\mu ^{\ast }\left( t\right) \right) \mathrm{d%
}W\left( t\right) , \\
-\mathrm{d}p_{2}\left( t\right) & = & \frac{\partial }{\partial x}\mathcal{H}%
_{2}\left( t,x\left( t\right) ,u_{2}\left( t\right) x\left( t\right)
+u_{1}\left( t\right) ,\mu ^{\ast }\left( t\right) ,p_{2}\left( t\right)
,q_{2}\left( t\right) \right) \mathrm{d}t-q_{2}\left( t\right) \mathrm{d}%
W\left( t\right) , \\
x\left( 0\right) & = & x_{0},\text{ }p_{2}\left( T\right) =\frac{\partial }{%
\partial x}\Phi _{2}\left( x\left( T\right) \right) ,%
\end{array}%
\right.  \label{acs}
\end{equation}%
with
\begin{eqnarray*}
\frac{\partial }{\partial x}\mathcal{H}_{2}\left( t,x,u_{2}x+u_{1},\mu
^{\ast },p_{2},q_{2}\right) &=&\left( \frac{\partial b}{\partial x}+\frac{%
\partial b}{\partial u}u_{2}\right) ^{\top }p_{2}+\left( \frac{\partial
\sigma }{\partial x}+\frac{\partial \sigma }{\partial u}u_{2}\right) ^{\top
}q_{2} \\
&&-\frac{\partial l_{2}}{\partial x}-u_{2}^{\top }\frac{\partial l_{2}}{%
\partial u}.
\end{eqnarray*}%
where $u_{2}$ and $u_{1}$ are adapted control processes with values in $%
\mathbb{R}^{m_{1}}$ and the ball $B_{K}(\mathbb{R}^{m_{1}\times n})$ with
radius $K$ in $\mathbb{R}^{m_{1}}$, respectively. We assume that the
coefficients $\frac{\partial }{\partial x}\mathcal{H}_{2}\left(
t,x,u_{2}x+u_{1},\mu ^{\ast },p_{2},q_{2}\right) ,$ $b$, $\sigma $ and $%
\frac{\partial }{\partial x}\Phi _{2}\left( x\right) $ satisfy the monotone
condition (A1)-(A2). Thus, the problem (\ref{necost})-(\ref{acs}) is
well-posed.

\begin{theorem}
\label{the2} Assume that the above coefficients of the problem \emph{(\ref%
{necost})-(\ref{acs})} satisfy the monotone conditions. Let $u^{\ast }\in
\mathcal{U}$ be an optimal solution to the leader's problem (\ref{cost2})-(%
\ref{a5}) with the corresponding state trajectory $\left( \bar{x}\left(
\cdot \right) ,\bar{p}_{2}\left( \cdot \right) ,\bar{q}_{2}\left( \cdot
\right) \right) \in \mathcal{N}^{2}\left[ 0,T\right] $, then there exists a
triple $\left( \chi \left( \cdot \right) ,p_{1}\left( \cdot \right)
,q_{1}\left( \cdot \right) \right) \in \mathcal{N}^{2}\left[ 0,T\right] $
such that%
\begin{eqnarray*}
&&\left( u^{\ast }\left( t,\bar{x}\left( t\right) \right) -\frac{\partial
u^{\ast }\left( t,\bar{x}\left( t\right) \right) }{\partial x}\bar{x}\left(
t\right) ,\frac{\partial u^{\ast }\left( t,\bar{x}\left( t\right) \right) }{%
\partial x}\right) \\
&=&\arg \min_{\left( u_{1},u_{2}\right) \in \mathbb{R}^{m_{1}}\times B_{K}(%
\mathbb{R}^{m_{1}\times n})}\mathcal{H}_{3}\left( t,u_{1},u_{2},\bar{x}%
\left( t\right) ,\chi \left( t\right) ,p_{1}\left( t\right) ,q_{1}\left(
t\right) ,\bar{p}_{2}\left( t\right) ,\bar{q}_{2}\left( t\right) \right)
\end{eqnarray*}%
subject to
\begin{equation}
\left\{
\begin{array}{rcl}
\mathrm{d}\chi \left( t\right) & = & -\frac{\partial }{\partial p_{2}}%
\mathcal{H}_{3}\mathrm{d}t-\frac{\partial }{\partial q_{2}}\mathcal{H}_{3}%
\mathrm{d}W\left( t\right) \\
\mathrm{d}p_{1}\left( t\right) & = & \frac{\partial }{\partial x}\mathcal{H}%
_{3}\mathrm{d}t+q_{1}\left( t\right) \mathrm{d}W\left( t\right) , \\
\chi \left( 0\right) & = & 0,\text{ }p_{1}\left( T\right) =-\frac{\partial
^{2}\Phi _{2}\left( \bar{x}\left( T\right) \right) }{\partial x^{2}}\chi
\left( T\right) +\frac{\partial \Phi _{1}\left( \bar{x}\left( T\right)
\right) }{\partial x},%
\end{array}%
\right.  \label{saclm}
\end{equation}%
where
\begin{eqnarray*}
&&\mathcal{H}_{3}\left( t,u_{1},u_{2},\bar{x},\chi ,p_{1},q_{1},\bar{p}_{2},%
\bar{q}_{2}\right) \\
&=&\left\langle p_{1},b\left( t,x,u_{2}x+u_{1}\left( t\right) ,v^{\ast
}\left( t,x,u_{2}x+u_{1},p_{2},q_{2}\right) \right) \right\rangle \\
&&+\left\langle q_{1},\sigma \left( t,x,u_{2}x+u_{1}\left( t\right) ,v^{\ast
}\left( t,x,u_{2}x+u_{1},p_{2},q_{2}\right) \right) \right\rangle \\
&&-\left\langle \chi ,\frac{\partial }{\partial x}\mathcal{H}_{2}\left(
t,x,u_{2}x+u_{1},\mu ^{\ast },p_{2},q_{2}\right) \right\rangle \\
&&+l_{1}\left( t,x,u_{2}x+u_{1},v^{\ast }\left(
t,x,u_{2}x+u_{1},p_{2},q_{2}\right) \right) .
\end{eqnarray*}%
Here $\frac{\partial }{\partial p_{2}}\mathcal{H}_{3},$ $\frac{\partial }{%
\partial q_{2}}\mathcal{H}_{3},$ $\frac{\partial }{\partial x}\mathcal{H}%
_{3} $ are evaluated at the point
\begin{equation*}
\left( t,u^{\ast }\left( t,\bar{x}\left( t\right) \right) -\frac{\partial
u^{\ast }\left( t,\bar{x}\left( t\right) \right) }{\partial x}\bar{x}\left(
t\right) ,\frac{\partial u^{\ast }\left( t,\bar{x}\left( t\right) \right) }{%
\partial x},\bar{x}\left( t\right) ,\chi \left( t\right) ,p_{1}\left(
t\right) ,q_{1}\left( t\right) ,\bar{p}_{2}\left( t\right) ,\bar{q}%
_{2}\left( t\right) \right) .
\end{equation*}
\end{theorem}

For reader's convenience, we present a brief proof as follows:

\paragraph{Proof.}

Obviously, from the definitions of the optimal costs $\mathcal{J}$ and $%
\mathcal{J}_{1}$ associated with problems (\ref{cost2}) and (\ref{necost}),
we have $\mathcal{J}_{1}\leq \mathcal{J}.$ Besides, based on $u^{\ast }$ of
problem (\ref{cost2}), we are able to structure a pair of control processes $%
\left( u_{1}^{\ast },u_{2}^{\ast }\right) $ for problem (\ref{necost}) as
follows:%
\begin{eqnarray}
u_{1}^{\ast }\left( t\right) &=&u^{\ast }\left( t,\bar{x}\left( t\right)
\right) -\frac{\partial u^{\ast }\left( t,\bar{x}\left( t\right) \right) }{%
\partial x}\bar{x}\left( t\right) ,  \label{ac1} \\
u_{2}^{\ast }\left( t\right) &=&\frac{\partial u^{\ast }\left( t,\bar{x}%
\left( t\right) \right) }{\partial x}.  \label{ac2}
\end{eqnarray}%
Employing these controls (\ref{ac1})-(\ref{ac2}), the FBSDEs (\ref{acs})
admit the same solution as that of (\ref{a5}) with the optimal strategy $%
u^{\ast }$ from which we derive
\begin{equation*}
u^{\ast }\left( t,\bar{x}\left( t\right) \right) =u_{2}^{\ast }\left(
t\right) \bar{x}\left( t\right) +u_{1}^{\ast }\left( t\right) .
\end{equation*}%
Consequently, $\mathcal{J}=\mathcal{J}_{1}$ and the above constructed $%
\left( u_{1}^{\ast },u_{2}^{\ast }\right) $ is an optimal control for
problem (\ref{necost}), leading to the same state trajectory $\left( \chi
\left( \cdot \right) ,p_{1}\left( \cdot \right) ,q_{1}\left( \cdot \right)
\right) .$ From the above arguments we state that%
\begin{equation*}
\left(
\begin{array}{cc}
u_{1}^{\ast }\left( t\right) =u^{\ast }\left( t,\bar{x}\left( t\right)
\right) -\frac{\partial u^{\ast }\left( t,\bar{x}\left( t\right) \right) }{%
\partial x}\bar{x}\left( t\right) , & u_{2}^{\ast }\left( t\right) =\frac{%
\partial u^{\ast }\left( t,\bar{x}\left( t\right) \right) }{\partial x}%
\end{array}%
\right)
\end{equation*}%
is indeed an optimal control for problem (\ref{necost}), whenever $u^{\ast }$
is a solution for problem (\ref{cost2}) with the corresponding forward state
$\bar{x}\left( t\right) .$ Therefore, we can establish the maximum principle
for problem (\ref{cost2}) of the leader via substituting (\ref{ac1})-(\ref%
{ac2}) into the necessary conditions satisfied by the optimal control for
problem (\ref{necost}).\hfill $\Box $

\begin{remark}
Note that whenever $u$ doesn't contain the state $x$, we claim that
Stackelberg solution is reduced to the AOL Stackelberg solution, and thus
the maximum principle for both cases coincides.
\end{remark}

\section{Application to linear quadratic Stackelberg games}

\label{sec3} In this section, the theoretical result obtained in Section \ref%
{sec2} will be applied to study linear quadratic Stackelberg games under the
AOL and ACLM information structures, respectively. Yong \cite{Yong2002}
derives stochastic Riccati equations for the follower and the leader
sequentially with random coefficients and diffusion term of the state
equation depending on controls, and the weight matrices in the cost
functionals are not necessarily positive definite. To be precise, the
follower gives his Riccati equation for any given strategy of the leader.
Then the leader solves his problem involving a system of FBSDEs, whose
coefficients depend on the solution of the follower's Riccati equation.
Finally, a further analysis of the state feedback representation of the
leader's optimal strategy provides the leader's Riccati equation. Under
certain conditions, the solvability of the leader's Riccati equation in the
case of deterministic coefficients is also discussed.

In contrast to Yong \cite{Yong2002}, we consider the similar system under
convex control constraint . The systems will contain the projection
operators which makes the system nonlinear (classical Riccati approach
fails). When supposing the control set is full space, we may let the
follower's Hamiltonian system as the leader's controlled state equation, and
consequently, the state feedback representation of the AOL Stackelberg
solution can be represented simultaneously for the leader and the follower.
As a result, the corresponding Riccati equation is different from that in
\cite{Yong2002}. Moreover, by means of a linear transformation to a standard
stochastic Riccati equation, we also prove that under certain conditions
there exists a unique solution to the Riccati equation with stochastic
coefficients studied by Tang \cite{Tang03}. For a linear quadratic
Stackelberg game under the ACLM case, we will see that the follower's
Hamiltonian system is no longer linear, and that keeps us from deriving a
Riccati equation if we handel the same way as in the AOL case. Instead, we
postulate that the forward variable $\chi $ is linear with respect to the
original state $x$, and then derive a kind of new but extremely complex
FBSDEs which plays the same role as the Riccati equation in the AOL case.

\subsection{The AOL fashion}

To this end, let us introduce the state equation and the cost functional for
leader and follower, respectively:
\begin{equation}
\left\{
\begin{array}{rcl}
\mathrm{d}x\left( t\right) & = & \left[ A\left( t\right) x\left( t\right)
+B_{1}\left( t\right) u\left( t\right) +B_{2}\left( t\right) v\left(
t\right) \right] \mathrm{d}t \\
&  & +\left[ C\left( t\right) x\left( t\right) +D_{1}\left( t\right) u\left(
t\right) +D_{2}\left( t\right) v\left( t\right) \right] \mathrm{d}W\left(
t\right) , \\
x\left( 0\right) & = & x_{0}\in \mathbb{R}^{n}.%
\end{array}%
\right.  \label{ex-sde1}
\end{equation}%
The cost functionals for the leader and the follower to minimize are given,
respectively, as follows:%
\begin{equation}
\mathcal{J}_{1}\left( u,v\right) =\frac{1}{2}\mathbb{E}\left[ \int_{0}^{T}%
\left[ \left( \left\langle Q_{1}\left( t\right) x\left( t\right) ,x\left(
t\right) \right\rangle +\left\langle R_{1}\left( t\right) u\left( t\right)
,u\left( t\right) \right\rangle \right) \right] \mathrm{d}t+\left\langle
\Phi _{1}x\left( T\right) ,x\left( T\right) \right\rangle \right]
\end{equation}%
and
\begin{equation}
\mathcal{J}_{2}\left( u,v\right) =\frac{1}{2}\mathbb{E}\left[
\int_{0}^{T}\left( \left\langle Q_{2}\left( t\right) x\left( t\right)
,x\left( t\right) \right\rangle +\left\langle R_{2}\left( t\right) v\left(
t\right) ,v\left( t\right) \right\rangle \right) \mathrm{d}t+\left\langle
\Phi _{2}x\left( T\right) ,x\left( T\right) \right\rangle \right] .
\end{equation}

We make the following three assumptions on the coefficients of the above
problem.

\begin{description}
\item[(H1)] Suppose that the matrix processes
\begin{eqnarray*}
A &:&\Omega \times \lbrack 0,T]\rightarrow \mathbb{R}^{n\times n},\text{ } \\
B_{1} &:&\Omega \times \lbrack 0,T]\rightarrow \mathbb{R}^{n\times m_{1}}, \\
B_{2} &:&\Omega \times \lbrack 0,T]\rightarrow \mathbb{R}^{n\times m_{2}}, \\
C &:&\Omega \times \lbrack 0,T]\rightarrow \mathbb{R}^{n\times n},\text{ } \\
D_{1} &:&\Omega \times \lbrack 0,T]\rightarrow \mathbb{R}^{n\times m_{1}},%
\text{ } \\
D_{2} &:&\Omega \times \lbrack 0,T]\rightarrow \mathbb{R}^{n\times m_{2}}, \\
Q_{1} &:&\Omega \times \lbrack 0,T]\rightarrow \mathbb{R}^{n\times n}, \\
\text{ }Q_{2} &:&\Omega \times \lbrack 0,T]\rightarrow \mathbb{R}^{n\times
n}, \\
R_{1} &:&\Omega \times \lbrack 0,T]\rightarrow \mathbb{R}^{m_{1}\times
m_{1}},\text{ } \\
R_{2} &:&\Omega \times \lbrack 0,T]\rightarrow \mathbb{R}^{m_{2}\times
m_{2}},
\end{eqnarray*}%
and the random matrices $\Phi _{1},\Phi _{2}:\Omega \rightarrow \mathbb{R}%
^{n}$ are uniformly bounded and $\left\{ \mathcal{F}_{t},0\leq t\leq
T\right\} $\}-adapted or $\mathcal{F}_{T}$ -measurable.

\item[(H2)] Suppose that the state weighting matrix process $Q_{1}$ and $%
Q_{2}$ are a.s. a.e. symmetric and nonnegative. Also suppose that the
terminal state weighting random matrix $\Phi _{1}$ and $\Phi _{2}$ are a.s.
symmetric and nonnegative.

\item[(H3)] Suppose that the control weighting matrix process $R_{1}$ and $%
R_{2}$ are a.s. a.e. symmetric and uniformly positive.
\end{description}

The Hamiltonian function can be expressed by
\begin{eqnarray}
\mathcal{H}_{2}\left( t,x,u,v,p_{2},q_{2}\right) &=&\left\langle
p_{2},A\left( t\right) x+B_{1}\left( t\right) u+B_{2}\left( t\right)
v\right\rangle  \notag \\
&&+\left\langle q_{2},C\left( t\right) x+D_{1}\left( t\right) u+D_{2}\left(
t\right) v\right\rangle  \notag \\
&&+\frac{1}{2}\left[ \left\langle Q_{2}\left( t\right) x,x\right\rangle
+\left\langle R_{2}\left( t\right) v,v\right\rangle \right] ,  \notag \\
\forall \left( t,x,u,v,p_{2},q_{2}\right) &\in &\left[ 0,T\right] \times
\mathbb{R}^{n}\times \mathbb{R}^{m_{1}}\times \mathbb{R}^{m_{2}}\times
\mathbb{R}^{n}\times \mathbb{R}^{n}.  \label{Hamiltonian function}
\end{eqnarray}%
The adjoint equation $\left( p_{2}\left( \cdot \right) ,q_{2}\left( \cdot
\right) \right) \in \mathcal{S}^{2}(0,T;\mathbb{R}^{n})\times \mathcal{M}%
^{2}(0,T;\mathbb{R}^{n\times d})$ such that%
\begin{equation}
\left\{
\begin{array}{rcl}
-\mathrm{d}p_{2}\left( t\right) & = & \left( A^{\top }\left( t\right)
p_{2}\left( t\right) +C^{\top }\left( t\right) q_{2}\left( t\right)
-Q_{2}\left( t\right) x\left( t\right) \right) \mathrm{d}t-q_{2}\left(
t\right) W\left( t\right) , \\
x\left( 0\right) & = & x_{0},\text{ }p_{2}\left( T\right) =-\Phi _{2}x\left(
T\right) .%
\end{array}%
\right.
\end{equation}%
Since $\Gamma _{2}$ is a closed convex set, then maximum principle reads as
the following local form
\begin{equation}
\left\langle -\frac{\partial \mathcal{H}_{2}}{\partial v}\left( t,x\left(
t\right) ,u\left( t\right) ,v^{\ast }\left( t\right) ,p_{2}\left( t\right)
,q_{2}\left( t\right) \right) ,v-v^{\ast }\left( t\right) \right\rangle \leq
0,\quad \forall v\in \Gamma _{2},\text{ a.e. }t\in \lbrack 0,T],\ \text{%
P-a.s.}  \label{convex maximum principle}
\end{equation}%
Hereafter, time argument is suppressed in case when no confusion occurs.

Noticing (\ref{Hamiltonian function}), then (\ref{convex maximum principle})
yields that
\begin{equation*}
\left\langle -B_{2}^{\top }p_{2}-D_{2}^{\top }q_{2}-R_{2}v^{\ast }\left(
t\right) ,v-v^{\ast }\left( t\right) \right\rangle \leq 0,\text{ for all }%
v\in \Gamma _{2},\text{ a.e. }t\in \lbrack 0,T],\ P\text{-a.s.}
\end{equation*}%
or equivalently (noticing $R_{2}>0$),
\begin{equation}
\left\langle R_{2}^{\frac{1}{2}}[-R_{2}^{-1}(B_{2}^{\top }p_{2}+D_{2}^{\top
}q_{2})-v^{\ast }\left( t\right) ],R_{2}^{\frac{1}{2}}(v-v^{\ast }\left(
t\right) )\right\rangle \leq 0,\text{ }\forall v\in \Gamma _{2},\text{ a.e. }%
t\in \lbrack 0,T],\ P\text{-a.s.}  \label{optimal control condition}
\end{equation}%
As $R_{2}(\cdot )>0$, we take the following norm on $\Gamma _{2}\subset
\mathbb{R}^{m_{2}}$ (which is equivalent to its Euclidean norm)
\begin{equation*}
\Vert x\Vert _{R_{2}}^{2}=\left\langle \left\langle x,x\right\rangle
\right\rangle :=\left\langle R_{2}^{\frac{1}{2}}x,R_{2}^{\frac{1}{2}%
}x\right\rangle ,
\end{equation*}%
and by the well-known results of convex analysis, we obtain that (\ref%
{optimal control condition}) is equivalent to
\begin{equation*}
v^{\ast }(t)=\mathbf{P}_{\Gamma _{2}}[-R_{2}^{-1}(t)(B_{2}^{\top
}(t)p_{2}(t)+D_{2}^{\top }(t)q_{2}(t))],\quad \text{ a.e. }t\in \lbrack
0,T],\ P\text{-a.s.},
\end{equation*}%
where $\mathbf{P}_{\Gamma _{2}}(\cdot )$ is the projection mapping from $%
\mathbb{R}^{m_{2}}$ to its closed convex subset $\Gamma _{2}$ under the norm
$\Vert \cdot \Vert _{R_{2}}$. For more details, see Appendix. From now on,
we denote
\begin{equation*}
\varphi _{2}(t,p,q):=\mathbf{P}_{\Gamma _{2}}[-R_{2}^{-1}(t)(B_{2}^{\top
}(t)p+D_{2}^{\top }(t)q)].
\end{equation*}%
The follower's optimal strategy as follows:%
\begin{equation*}
v^{\ast }\left( t\right) =\varphi _{2}(t,p_{1}\left( t\right) ,q_{1}\left(
t\right) ):=\mathbf{P}_{\Gamma _{2}}[-R_{2}^{-1}(t)(B_{2}^{\top
}(t)p_{2}\left( t\right) +D_{2}^{\top }(t)q_{2}\left( t\right) )].
\end{equation*}%
%
%
%
%
%
%
%
%
%
%
%
%
%
%
%
%
%
%
%
%
%
%
%
%
%
%
%
%
%
%
%
%
%
%
%
%
%
%
%
%
%
%
%
%
%
%
%
%
%
%
%
%
%
%
%
%
%
%
%
%
%
%
%
%
%
%
%
%
%
%
%
%
%
%
%
%
%
%
%
%
%
%
%
%
%
%
%
%
%
%
%
%
%
%
%
%
%
%
%
%
%
%
%
%
%
%
%
%
%
%
%
%
%
%
%
%
%
%
%
%
%
%
%
%
%
%
%
%
%
%
%
%
%
%
%

Now, we focus on the leader's problem. Her/His aim is to seek an optimal
control $u^{\ast }\left( \cdot \right) $ such that
\begin{equation*}
\mathcal{J}_{1}\left( u^{\ast },v^{\ast }\right) =\inf_{u\left( \cdot
\right) \in \mathcal{U}}\mathcal{J}_{1}\left( u,v^{\ast }\right)
\end{equation*}%
subject to

\begin{equation}
\left\{
\begin{array}{rcl}
\mathrm{d}x\left( t\right) & = & \left[ A\left( t\right) x\left( t\right)
+B_{1}\left( t\right) u\left( t\right) +B_{2}\left( t\right) \varphi
_{2}(t,p_{2}\left( t\right) ,q_{2}\left( t\right) )\right] \mathrm{d}t \\
&  & +\left[ C\left( t\right) x\left( t\right) +D_{1}\left( t\right) u\left(
t\right) +D_{2}\left( t\right) \varphi _{2}(t,p_{2}\left( t\right)
,q_{2}\left( t\right) )\right] \mathrm{d}W\left( t\right) , \\
-\mathrm{d}p_{2}\left( t\right) & = & \left[ A^{\top }\left( t\right)
p_{2}\left( t\right) +C^{\top }\left( t\right) q_{2}\left( t\right)
+Q_{2}\left( t\right) x\left( t\right) \right] \mathrm{d}t-q_{2}\left(
t\right) W\left( t\right) , \\
x\left( 0\right) & = & x_{0},\text{ }p_{2}\left( T\right) =\Phi _{2}x\left(
T\right) .%
\end{array}%
\right.  \label{state1}
\end{equation}%
Obviously, FBSDEs (\ref{state1}) are fully coupled, which contains a
nonlinear term $\varphi _{2}(t,p_{2}\left( t\right) ,q_{2}\left( t\right) ).$
Nonetheless, under certain assumptions, we are able to prove the existence
and uniqueness of such equations.

\begin{theorem}
\label{the1} Assume that \emph{(H1)-(H3)} are in force. Then, for any given $%
u\left( \cdot \right) \in \mathcal{U},$ FBSDEs (\ref{state1}) admit a unique
adapted solution $\left( x\left( \cdot \right) ,p_{2}\left( \cdot \right)
,q_{2}\left( \cdot \right) \right) \in \mathcal{N}^{2}\left[ 0,T\right] .$
\end{theorem}

The proof can be found in the Appendix \ref{APP}.

\begin{remark}
Due to the nonlinearity of (\ref{state1}), the classical approach of Riccati
equation is not applicable in this case. Moreover, the methodology developed
in Hu and Zhou \cite{hz} can not applied directly. On the one hand, the
control domain there is postulated to be a closed cone involving the
original point. From (5.2) in \cite{hz}, we know that the optimal feedback
control can be expressed as control process multiplying by the state
variable. On the other hand, note that the liner system (\ref{ex-sde1}) is
non-homogeneous linear equation, which doesn't satisfy the framework in \cite%
{hz} since the equations (5.24) and (5.25) in Hu and Zhou \cite{hz} can be
represented explicitly. However, in our paper, we have two controls
simultaneously, whose system, of course, is \textit{non-homogeneous.}
\end{remark}

Now we are ready to find the optimal control for leader. The leader's
problem is well-posed since for every $u\left( \cdot \right) \in \mathcal{U}$%
, the FBSDEs (\ref{state1}) admits a unique solution. From Proposition \ref%
{pro1}, it is easy to derive the leader's optimal strategy as follows:%
\begin{equation*}
u^{\ast }\left( t\right) =\varphi _{1}(t,p_{1}\left( t\right) ,q_{1}\left(
t\right) ):=\mathbf{P}_{\Gamma _{1}}[-R_{1}^{-1}(t)(B_{1}^{\top
}(t)p_{1}\left( t\right) +D_{1}^{\top }(t)q_{1}\left( t\right) )],
\end{equation*}%
where
\begin{equation}
\left\{
\begin{array}{rcl}
\mathrm{d}k\left( t\right) & = & -\frac{\partial }{\partial p_{2}}\mathcal{H}%
_{1}\mathrm{d}t-\frac{\partial }{\partial q_{2}}\mathcal{H}_{1}\mathrm{d}%
W\left( t\right) \\
\mathrm{d}p_{1}\left( t\right) & = & -\frac{\partial }{\partial x}\mathcal{H}%
_{1}\mathrm{d}t+q_{1}\left( t\right) \mathrm{d}W\left( t\right) , \\
k\left( 0\right) & = & 0,\text{ }p_{1}\left( T\right) =-\frac{\partial ^{2}}{%
\partial x^{2}}\Phi _{2}\left( x\left( T\right) \right) k\left( T\right) +%
\frac{\partial }{\partial x}\Phi _{1}\left( x\left( T\right) \right) ,%
\end{array}%
\right.
\end{equation}%
and
\begin{eqnarray*}
&&\mathcal{H}_{1}\left( t,u,x,k,p_{1},p_{2},q_{1},q_{2}\right) \\
&=&\left\langle p_{1},Ax+B_{1}u+B_{2}\varphi _{2}(t,p_{2},q_{2})\right\rangle
\\
&&+\left\langle q_{1},Cx+D_{1}u+D_{2}\varphi _{2}(t,p_{2},q_{2})\right\rangle
\\
&&+\frac{1}{2}\left( \left\langle Q_{1}x,x\right\rangle +\left\langle
R_{1}u,u\right\rangle \right) \\
&&-\left\langle k,A^{\top }p_{2}+C^{\top }q_{2}+Q_{2}x\right\rangle .
\end{eqnarray*}%
From the uniqueness of the optimal strategy and Proposition \ref{projection
theorem}, we also know that the FBSDEs:%
\begin{equation}
\left\{
\begin{array}{rcl}
\mathrm{d}x\left( t\right) & = & \left[ A\left( t\right) x\left( t\right)
+B_{1}\left( t\right) u\left( t\right) +B_{2}\left( t\right) \varphi
_{2}(t,p_{2}\left( t\right) ,q_{2}\left( t\right) )\right] \mathrm{d}t \\
&  & +\left[ C\left( t\right) x\left( t\right) +D_{1}\left( t\right) u\left(
t\right) +D_{2}\left( t\right) \varphi _{2}(t,p_{2}\left( t\right)
,q_{2}\left( t\right) )\right] \mathrm{d}W\left( t\right) , \\
-\mathrm{d}p_{2}\left( t\right) & = & \left[ A^{\top }\left( t\right)
p_{2}\left( t\right) +C^{\top }\left( t\right) q_{2}\left( t\right)
+Q_{2}\left( t\right) x\left( t\right) \right] \mathrm{d}t-q_{2}\left(
t\right) W\left( t\right) , \\
\mathrm{d}k\left( t\right) & = & \left[ -B_{2}\left( t\right) \frac{\partial
\varphi _{2}(t,p_{2}\left( t\right) ,q_{2}\left( t\right) )}{\partial p_{2}}%
p_{1}\left( t\right) -D_{2}\left( t\right) \frac{\partial \varphi
_{2}(t,p_{2}\left( t\right) ,q_{2}\left( t\right) )}{\partial p_{2}}%
q_{1}\left( t\right) +A\left( t\right) k\left( t\right) \right] \mathrm{d}t
\\
&  & +\left[ -B_{2}\left( t\right) \frac{\partial \varphi _{2}(t,p_{2}\left(
t\right) ,q_{2}\left( t\right) )}{\partial q_{2}}p_{1}\left( t\right)
-D_{2}\left( t\right) \frac{\partial \varphi _{2}(t,p_{2}\left( t\right)
,q_{2}\left( t\right) )}{\partial q_{2}}q_{1}\left( t\right) +C\left(
t\right) k\left( t\right) \right] \mathrm{d}W\left( t\right) , \\
-\mathrm{d}p_{1}\left( t\right) & = & \left[ A^{\top }\left( t\right)
p_{1}\left( t\right) +C^{\top }\left( t\right) q_{1}\left( t\right)
+Q_{1}\left( t\right) x\left( t\right) -Q_{2}\left( t\right) k\left(
t\right) \right] \mathrm{d}t-q_{1}\left( t\right) \mathrm{d}W\left( t\right)
, \\
x\left( 0\right) & = & x_{0},\text{ }k\left( 0\right) =0,\text{ }p_{1}\left(
T\right) =\Phi _{2}k\left( T\right) -\Phi _{1}x\left( T\right) ,\text{ }%
p_{2}\left( T\right) =-\Phi _{2}x\left( T\right)%
\end{array}%
\right.  \label{fc}
\end{equation}%
has a unique solution under the assumption $\frac{\partial }{\partial p_{2}}%
\varphi _{2}(t,p_{2}\left( t\right) ,q_{2}\left( t\right) )$ and $\frac{%
\partial }{\partial q_{2}}\varphi _{2}(t,p_{2}\left( t\right) ,q_{2}\left(
t\right) )$ are well-defined.

Finally, we have the following coupled systems:
\begin{equation}
\left\{
\begin{array}{rcl}
\mathrm{d}x\left( t\right) & = & \left[ A\left( t\right) x\left( t\right)
+B_{1}\left( t\right) \varphi _{1}(t,p_{1}\left( t\right) ,q_{1}\left(
t\right) )+B_{2}\left( t\right) \varphi _{2}(t,p_{2}\left( t\right)
,q_{2}\left( t\right) )\right] \mathrm{d}t \\
&  & +\left[ C\left( t\right) x\left( t\right) +D_{1}\left( t\right) \varphi
_{1}(t,p_{1}\left( t\right) ,q_{1}\left( t\right) )+D_{2}\left( t\right)
\varphi _{2}(t,p_{2}\left( t\right) ,q_{2}\left( t\right) )\right] \mathrm{d}%
W\left( t\right) , \\
-\mathrm{d}p_{2}\left( t\right) & = & \left[ A^{\top }\left( t\right)
p_{2}\left( t\right) +C^{\top }\left( t\right) q_{2}\left( t\right)
+Q_{2}\left( t\right) x\left( t\right) \right] \mathrm{d}t-q_{2}\left(
t\right) W\left( t\right) , \\
\mathrm{d}k\left( t\right) & = & \left[ -B_{2}\left( t\right) \frac{\partial
\varphi _{2}(t,p_{2}\left( t\right) ,q_{2}\left( t\right) )}{\partial p_{2}}%
p_{1}\left( t\right) -D_{2}\left( t\right) \frac{\partial \varphi
_{2}(t,p_{2}\left( t\right) ,q_{2}\left( t\right) )}{\partial p_{2}}%
q_{1}\left( t\right) +A\left( t\right) k\left( t\right) \right] \mathrm{d}t
\\
&  & +\left[ -B_{2}\left( t\right) \frac{\partial \varphi _{2}(t,p_{2}\left(
t\right) ,q_{2}\left( t\right) )}{\partial q_{2}}p_{1}\left( t\right)
-D_{2}\left( t\right) \frac{\partial \varphi _{2}(t,p_{2}\left( t\right)
,q_{2}\left( t\right) )}{\partial q_{2}}q_{1}\left( t\right) +C\left(
t\right) k\left( t\right) \right] \mathrm{d}W\left( t\right) , \\
-\mathrm{d}p_{1}\left( t\right) & = & \left[ A^{\top }\left( t\right)
p_{1}\left( t\right) +C^{\top }\left( t\right) q_{1}\left( t\right)
+Q_{1}\left( t\right) x\left( t\right) -Q_{2}\left( t\right) k\left(
t\right) \right] \mathrm{d}t-q_{1}\left( t\right) \mathrm{d}W\left( t\right)
, \\
x\left( 0\right) & = & x_{0},\text{ }k\left( 0\right) =0,\text{ }p_{1}\left(
T\right) =\Phi _{2}k\left( T\right) -\Phi _{1}x\left( T\right) ,\text{ }%
p_{2}\left( T\right) =-\Phi _{2}x\left( T\right)%
\end{array}%
\right.  \label{fc2}
\end{equation}

Next we set $\Gamma _{1}=\mathbb{R}^{m_{1}},\Gamma _{2}=\mathbb{R}^{m_{2}}.$
We observe that the Riccati equation approach is really applicable in this
case, and the AOL Stackelberg solution $\left( u^{\ast },v^{\ast }\right) $
can be written as%
\begin{equation*}
\left\{
\begin{array}{l}
u^{\ast }\left( t\right) =-R_{1}^{-1}(t)(B_{1}^{\top }(t)p_{1}\left(
t\right) +D_{1}^{\top }(t)q_{1}\left( t\right) ), \\
v^{\ast }\left( t\right) =-R_{2}^{-1}(t)(B_{2}^{\top }(t)p_{2}\left(
t\right) +D_{2}^{\top }(t)q_{2}\left( t\right) ).%
\end{array}%
\right.
\end{equation*}%
In this case, FBSDEs (\ref{fc2}) turn into
\begin{equation}
\left\{
\begin{array}{rcl}
\mathrm{d}x\left( t\right) & = & \left[ A\left( t\right) x\left( t\right)
+B_{1}\left( t\right) \varphi _{1}(t,p_{1}\left( t\right) ,q_{1}\left(
t\right) )+B_{2}\left( t\right) \varphi _{2}(t,p_{2}\left( t\right)
,q_{2}\left( t\right) )\right] \mathrm{d}t \\
&  & +\left[ C\left( t\right) x\left( t\right) +D_{1}\left( t\right) \varphi
_{1}(t,p_{1}\left( t\right) ,q_{1}\left( t\right) )+D_{2}\left( t\right)
\varphi _{2}(t,p_{2}\left( t\right) ,q_{2}\left( t\right) )\right] \mathrm{d}%
W\left( t\right) , \\
\mathrm{d}k\left( t\right) & = & \left[ B_{2}\left( t\right)
R_{2}^{-1}(t)B_{2}^{\top }(t)p_{1}\left( t\right) +D_{2}\left( t\right)
R_{2}^{-1}(t)B_{2}^{\top }(t)q_{1}\left( t\right) +A\left( t\right) k\left(
t\right) \right] \mathrm{d}t \\
&  & +\left[ B_{2}\left( t\right) R_{2}^{-1}(t)D_{2}^{\top }(t)p_{1}\left(
t\right) +D_{2}\left( t\right) R_{2}^{-1}(t)D_{2}^{\top }(t)q_{1}\left(
t\right) +C\left( t\right) k\left( t\right) \right] \mathrm{d}W\left(
t\right) , \\
-\mathrm{d}p_{1}\left( t\right) & = & \left[ A^{\top }\left( t\right)
p_{1}\left( t\right) +C^{\top }\left( t\right) q_{1}\left( t\right)
+Q_{1}\left( t\right) x\left( t\right) -Q_{2}\left( t\right) k\left(
t\right) \right] \mathrm{d}t-q_{1}\left( t\right) \mathrm{d}W\left( t\right)
, \\
-\mathrm{d}p_{2}\left( t\right) & = & \left[ A^{\top }\left( t\right)
p_{2}\left( t\right) +C^{\top }\left( t\right) q_{2}\left( t\right)
+Q_{2}\left( t\right) x\left( t\right) \right] \mathrm{d}t-q_{2}\left(
t\right) W\left( t\right) , \\
x\left( 0\right) & = & x_{0},\text{ }k\left( 0\right) =0,\text{ }p_{1}\left(
T\right) =-\Phi _{2}k\left( T\right) +\Phi _{1}x\left( T\right) ,\text{ }%
p_{2}\left( T\right) =\Phi _{2}x\left( T\right) .%
\end{array}%
\right.
\end{equation}%
It is possible to derive the feedback representation of the Stackelberg
solution $\left( u^{\ast },v^{\ast }\right) $ in terms of the state $\left(
x,k\right) $. To this end, we introduce the following notations for
simplicity,%
\begin{equation*}
X=\left(
\begin{array}{c}
x \\
k%
\end{array}%
\right) ,\text{ }P=\left(
\begin{array}{c}
p_{1} \\
p_{2}%
\end{array}%
\right) ,\text{ }Q=\left(
\begin{array}{c}
q_{1} \\
q_{2}%
\end{array}%
\right)
\end{equation*}%
and
\begin{eqnarray*}
\mathcal{A} &=&\left(
\begin{array}{cc}
A & 0 \\
0 & A%
\end{array}%
\right) ,\text{ }\mathcal{C}=\left(
\begin{array}{cc}
C & 0 \\
0 & C%
\end{array}%
\right) ,\text{ }X_{0}=\left(
\begin{array}{c}
x_{0} \\
0%
\end{array}%
\right) , \\
\mathcal{B}_{1} &=&\left(
\begin{array}{cc}
B_{1}R_{1}^{-1}B_{1}^{\top } & B_{2}R_{2}^{-1}B_{2}^{\top } \\
-B_{2}R_{2}^{-1}B_{2}^{\top } & 0%
\end{array}%
\right) ,\text{ }\mathcal{B}_{2}=\left(
\begin{array}{cc}
B_{1}R_{1}^{-1}D_{1}^{\top } & B_{2}R_{2}^{-1}D_{2}^{\top } \\
-D_{2}R_{2}^{-1}B_{2}^{\top } & 0%
\end{array}%
\right) , \\
\mathcal{D}_{1} &=&\left(
\begin{array}{cc}
D_{1}R_{1}^{-1}B_{1}^{\top } & D_{2}R_{2}^{-1}B_{2}^{\top } \\
-B_{2}R_{2}^{-1}D_{2}^{\top } & 0%
\end{array}%
\right) ,\text{ }\mathcal{D}_{2}=\left(
\begin{array}{cc}
D_{1}R_{1}^{-1}D_{1}^{\top } & D_{2}R_{2}^{-1}D_{2}^{\top } \\
-D_{2}R_{2}^{-1}D_{2}^{\top } & 0%
\end{array}%
\right) , \\
Q_{1} &=&\left(
\begin{array}{cc}
Q_{1} & -Q_{2} \\
Q_{2} & 0%
\end{array}%
\right) ,\text{ }\Phi _{1}=\left(
\begin{array}{cc}
\Phi _{1} & -\Phi _{2} \\
\Phi _{2} & 0%
\end{array}%
\right) .
\end{eqnarray*}%
Then, FBSDEs (\ref{fc2}) can be rewritten as
\begin{equation}
\left\{
\begin{array}{rcl}
\mathrm{d}X\left( t\right) & = & \left[ \mathcal{A}\left( t\right) X\left(
t\right) -\mathcal{B}_{1}\left( t\right) P\left( t\right) -\mathcal{B}%
_{2}\left( t\right) Q\left( t\right) \right] \mathrm{d}t \\
&  & +\left[ \mathcal{C}\left( t\right) X\left( t\right) -\mathcal{D}%
_{1}\left( t\right) P\left( t\right) -\mathcal{D}_{2}\left( t\right) Q\left(
t\right) \right] \mathrm{d}W\left( t\right) , \\
-\mathrm{d}P\left( t\right) & = & \left[ \mathcal{A}\left( t\right) P\left(
t\right) +\mathcal{C}^{\top }\left( t\right) Q\left( t\right) +Q_{1}\left(
t\right) X\left( t\right) \right] \mathrm{d}t-Q\left( t\right) \mathrm{d}%
W\left( t\right) , \\
X\left( 0\right) & = & X_{0},\text{ }P\left( T\right) =\hat{\Phi}X\left(
T\right) .%
\end{array}%
\right.  \label{sfbsde}
\end{equation}%
We are ready to derive the Riccati equation. To this end, assume that there
exists a matrix-valued process $\mathcal{P}$ such that
\begin{equation*}
P\left( t\right) =\mathcal{R}\left( t\right) X\left( t\right) ,\text{ }t\in %
\left[ 0,T\right] ,
\end{equation*}%
where $\mathcal{R}\left( t\right) $ is an $\mathcal{F}_{t}$-adapted process
with values in $\mathbb{R}^{n\times n}$. In general, $\mathcal{R}\left(
t\right) $ is not a bounded variation function with respect to $t$. We
tentatively assume that $\mathcal{R}\left( t\right) $ is a semi-martingale
\begin{equation}
\mathcal{R}\left( t\right) =\Phi _{1}+\int_{t}^{T}\Pi \left( s\right)
\mathrm{d}s-\int_{t}^{T}\Psi \left( s\right) \mathrm{d}W\left( s\right) ,%
\text{ }0\leq t\leq T,  \label{exp}
\end{equation}%
Applying the It\^{o}'s formula to $\mathcal{R}\left( \cdot \right) X\left(
\cdot \right) ,$ we obtain
\begin{eqnarray}
&&\mathcal{R}\left( t\right) \left[ \mathcal{A}\left( t\right) X\left(
t\right) -\mathcal{B}_{1}\left( t\right) \mathcal{R}\left( t\right) X\left(
t\right) -\mathcal{B}_{2}\left( t\right) Q\left( t\right) \right] \mathrm{d}t
\notag \\
&&-\Pi \left( t\right) X\left( t\right) \mathrm{d}t+\Psi \left( t\right) %
\left[ \mathcal{C}\left( t\right) X\left( t\right) -\mathcal{D}_{1}\left(
t\right) \mathcal{R}\left( t\right) X\left( t\right) -\mathcal{D}_{2}\left(
t\right) Q\left( t\right) \right] \mathrm{d}t  \notag \\
&&+\mathcal{R}\left( t\right) \left[ \mathcal{C}\left( t\right) X\left(
t\right) -\mathcal{D}_{1}\left( t\right) \mathcal{R}\left( t\right) X\left(
t\right) -\mathcal{D}_{2}\left( t\right) Q\left( t\right) \right] \mathrm{d}%
W\left( t\right) +\Psi \left( t\right) X\left( t\right) \mathrm{d}W\left(
t\right)  \notag \\
&=&\mathrm{d}P\left( t\right)  \notag \\
&=&-\left[ \mathcal{A}\left( t\right) \mathcal{R}\left( t\right) X\left(
t\right) +\mathcal{C}^{\top }\left( t\right) Q\left( t\right) +Q_{1}\left(
t\right) X\left( t\right) \right] \mathrm{d}t+Q\left( t\right) \mathrm{d}%
W\left( t\right) .  \label{ricar}
\end{eqnarray}%
It is easy to see
\begin{equation*}
Q\left( t\right) =\mathcal{R}\left( t\right) \left[ \mathcal{C}\left(
t\right) X\left( t\right) -\mathcal{D}_{1}\left( t\right) \mathcal{R}\left(
t\right) X\left( t\right) -\mathcal{D}_{2}\left( t\right) Q\left( t\right) %
\right] +\Psi \left( t\right) X\left( t\right) ,
\end{equation*}%
from which we get%
\begin{equation}
Q\left( t\right) =\Xi \left( t\right) X\left( t\right) ,  \label{inser}
\end{equation}%
where
\begin{equation*}
\Xi \left( t\right) =\left( I+\mathcal{R}\left( t\right) \mathcal{D}%
_{2}\left( t\right) \right) ^{-1}\left[ \mathcal{R}\left( t\right) \mathcal{C%
}\left( t\right) -\mathcal{R}\left( t\right) \mathcal{D}_{1}\left( t\right)
\mathcal{R}\left( t\right) +\Psi \left( t\right) \right] .
\end{equation*}%
Inserting (\ref{inser}) into (\ref{ricar}), we have
\begin{eqnarray*}
&&\mathcal{R}\left( t\right) \left[ \mathcal{A}\left( t\right) -\mathcal{B}%
_{1}\left( t\right) \mathcal{R}\left( t\right) -\mathcal{B}_{2}\left(
t\right) \Xi \left( t\right) \right] X\left( t\right) \mathrm{d}t \\
&&-\Pi \left( t\right) X\left( t\right) \mathrm{d}t+\Psi \left( t\right) %
\left[ \mathcal{C}\left( t\right) -\mathcal{D}_{1}\left( t\right) \mathcal{R}%
\left( t\right) -\mathcal{D}_{2}\left( t\right) \Xi \left( t\right) \right]
X\left( t\right) \mathrm{d}t \\
&=&-\left[ \mathcal{A}\left( t\right) \mathcal{R}\left( t\right) +\mathcal{C}%
^{\top }\left( t\right) \Xi \left( t\right) +Q_{1}\left( t\right) \right]
X\left( t\right) \mathrm{d}t,
\end{eqnarray*}%
which yields%
\begin{eqnarray*}
\Pi \left( t\right) &=&\mathcal{A}\left( t\right) \mathcal{R}\left( t\right)
+\mathcal{C}^{\top }\left( t\right) \Xi \left( t\right) +Q_{1}\left( t\right)
\\
&&+\mathcal{R}\left( t\right) \left[ \mathcal{A}\left( t\right) -\mathcal{B}%
_{1}\left( t\right) \mathcal{R}\left( t\right) -\mathcal{B}_{2}\left(
t\right) \Xi \left( t\right) \right] \\
&&+\Psi \left( t\right) \left[ \mathcal{C}\left( t\right) -\mathcal{D}%
_{1}\left( t\right) \mathcal{R}\left( t\right) -\mathcal{D}_{2}\left(
t\right) \Xi \left( t\right) \right] .
\end{eqnarray*}%
Consequently, we obtain
\begin{eqnarray}
\mathrm{d}\mathcal{R}\left( t\right) &=&[\mathcal{A}\left( t\right) \mathcal{%
R}\left( t\right) +\mathcal{C}^{\top }\left( t\right) \Xi \left( t\right)
+Q_{1}\left( t\right)  \notag \\
&&+\mathcal{R}\left( t\right) \left[ \mathcal{A}\left( t\right) -\mathcal{B}%
_{1}\left( t\right) \mathcal{R}\left( t\right) -\mathcal{B}_{2}\left(
t\right) \Xi \left( t\right) \right]  \notag \\
&&+\Psi \left( t\right) \left[ \mathcal{C}\left( t\right) -\mathcal{D}%
_{1}\left( t\right) \mathcal{R}\left( t\right) -\mathcal{D}_{2}\left(
t\right) \Xi \left( t\right) \right] ]\mathrm{d}t  \notag \\
&&-\Psi \left( t\right) \mathrm{d}W\left( t\right)  \notag \\
\mathcal{R}\left( T\right) &=&\Phi _{1}.  \label{sri}
\end{eqnarray}

\begin{remark}
Whenever $D_{1}\left( t\right) =D_{2}\left( t\right) =0,$ the Riccati
equation (\ref{sri}) becomes the right form (5.21) in Bensoussan et al. \cite%
{BCS2015}.
\end{remark}

We should claim that the Riccati equation (\ref{sri}) is just another
equivalent form in Tang \cite{Tang03} (see Discussion in Appendix).
Nonetheless, the coefficients $\mathcal{B}_{1},$ $\mathcal{B}_{2},$ $%
\mathcal{D}_{1},$ $\mathcal{D}_{2}$ and $Q_{1}$ are not symmetric matrices.
Next, we shall introduce a linear transformation to turn (\ref{sri}) into a
standard Riccati equation for $n=1.$

\begin{theorem}
For $n=1,$ we assume that
\begin{eqnarray*}
\frac{Q_{2}}{Q_{1}} &=&\frac{\Phi _{2}}{\Phi _{1}}, \\
\frac{B_{2}R_{2}^{-1}B_{2}^{\top }}{B_{1}R_{1}^{-1}B_{1}^{\top }} &=&\frac{%
B_{2}R_{2}^{-1}D_{2}^{\top }}{B_{1}R_{1}^{-1}D_{1}^{\top }}=\frac{%
D_{2}R_{2}^{-1}B_{2}^{\top }}{D_{1}R_{1}^{-1}B_{1}^{\top }}=\frac{%
D_{2}R_{2}^{-1}D_{2}^{\top }}{D_{1}R_{1}^{-1}D_{1}^{\top }}
\end{eqnarray*}
are in force. Then the Riccati equation (\ref{sri}) admits a unique solution.
\end{theorem}

\paragraph{Proof.}

Let
\begin{eqnarray}
\frac{Q_{2}}{Q_{1}} &=&\frac{\Phi _{2}}{\Phi _{1}}=\lambda ,\text{ }  \notag
\\
\frac{B_{2}R_{2}^{-1}B_{2}^{\top }}{B_{1}R_{1}^{-1}B_{1}^{\top }} &=&\frac{%
B_{2}R_{2}^{-1}D_{2}^{\top }}{B_{1}R_{1}^{-1}D_{1}^{\top }}=\frac{%
D_{2}R_{2}^{-1}B_{2}^{\top }}{D_{1}R_{1}^{-1}B_{1}^{\top }}=\frac{%
D_{2}R_{2}^{-1}D_{2}^{\top }}{D_{1}R_{1}^{-1}D_{1}^{\top }}=\mu .
\label{syass}
\end{eqnarray}%
Let us introduce a linear transformation
\begin{equation*}
X=\bar{X},P=\Upsilon \bar{P},\text{ }Q=\Upsilon \bar{Q},
\end{equation*}%
via a matrix $\Upsilon $ (determined later). So The FBSDEs (\ref{sfbsde})
can be expressed as
\begin{equation}
\left\{
\begin{array}{rcl}
\mathrm{d}\bar{X}\left( t\right) & = & \left[ \mathcal{\bar{A}}\left(
t\right) \bar{X}\left( t\right) -\mathcal{\bar{B}}_{1}\left( t\right) \bar{P}%
\left( t\right) -\mathcal{\bar{B}}_{2}\left( t\right) \bar{Q}\left( t\right) %
\right] \mathrm{d}t \\
&  & +\left[ \mathcal{\bar{C}}\left( t\right) \bar{X}\left( t\right) -%
\mathcal{\bar{D}}_{1}\left( t\right) \bar{P}\left( t\right) -\mathcal{\bar{D}%
}_{2}\left( t\right) \bar{Q}\left( t\right) \right] \mathrm{d}W\left(
t\right) , \\
-\mathrm{d}\bar{P}\left( t\right) & = & \left[ \mathcal{\bar{A}}\left(
t\right) \bar{P}\left( t\right) +\mathcal{\bar{C}}^{\top }\left( t\right)
\bar{Q}\left( t\right) +\bar{Q}_{1}\left( t\right) \bar{X}\left( t\right) %
\right] \mathrm{d}t-Q\left( t\right) \mathrm{d}W\left( t\right) , \\
\bar{X}\left( 0\right) & = & X_{0},\text{ }\bar{P}\left( T\right) =\bar{\Phi}%
_{1}\bar{X}\left( T\right) ,%
\end{array}%
\right.
\end{equation}%
where
\begin{eqnarray*}
\mathcal{\bar{A}} &=&\Upsilon ^{-1}A\Upsilon ,\text{ }\mathcal{\bar{B}}_{1}=%
\mathcal{B}_{1}\Upsilon ,\text{ }\mathcal{\bar{B}}_{2}=\mathcal{B}%
_{2}\Upsilon , \\
\mathcal{\bar{C}} &=&\Upsilon ^{-1}C\Upsilon ,\text{ }\mathcal{\bar{D}}_{1}=%
\mathcal{D}_{1}\Upsilon ,\text{ }\mathcal{\bar{D}}_{2}=\mathcal{D}%
_{2}\Upsilon , \\
\bar{Q}_{1} &=&\Upsilon ^{-1}Q_{1},\text{ }\bar{\Phi}_{1}=\Upsilon ^{-1}\Phi
_{1}.
\end{eqnarray*}%
Note that $A,C$ are symmetric metrics. We are going to seek $\Upsilon $ such
that $\mathcal{\bar{B}}_{1},\mathcal{\bar{B}}_{2},\mathcal{\bar{D}}_{1},%
\mathcal{\bar{D}}_{2},\bar{Q}_{1}$ and $\bar{\Phi}_{1}$ are symmetric. From
the assumption (\ref{syass}), it is easy to compute that
\begin{equation*}
\Upsilon =\left(
\begin{array}{cc}
1 & -2\mu \\
2\lambda & 1%
\end{array}%
\right) .
\end{equation*}%
Therefore, $\mathcal{\bar{A}=A}$, $\mathcal{\bar{C}=C}$,
\begin{eqnarray*}
\mathcal{\bar{B}}_{1} &=&\left(
\begin{array}{cc}
B_{1}R_{1}^{-1}B_{1}^{\top }+2\lambda B_{2}R_{2}^{-1}B_{2}^{\top } &
-B_{2}R_{2}^{-1}B_{2}^{\top } \\
-B_{2}R_{2}^{-1}B_{2}^{\top } & 2\mu B_{2}R_{2}^{-1}B_{2}^{\top }%
\end{array}%
\right) , \\
\mathcal{\bar{B}}_{2} &=&\left(
\begin{array}{cc}
B_{1}R_{1}^{-1}D_{1}^{\top }+2\lambda B_{2}R_{2}^{-1}D_{2}^{\top } &
-B_{2}R_{2}^{-1}B_{2}^{\top } \\
-D_{2}R_{2}^{-1}B_{2}^{\top } & 2\mu D_{2}R_{2}^{-1}B_{2}^{\top }%
\end{array}%
\right) , \\
\mathcal{\bar{D}}_{1} &=&\left(
\begin{array}{cc}
D_{1}R_{1}^{-1}B_{1}^{\top }+2\lambda D_{2}R_{2}^{-1}B_{2}^{\top } &
-B_{2}R_{2}^{-1}D_{2}^{\top } \\
-B_{2}R_{2}^{-1}D_{2}^{\top } & 2\mu B_{2}R_{2}^{-1}D_{2}^{\top }%
\end{array}%
\right) , \\
\mathcal{\bar{D}}_{2} &=&\left(
\begin{array}{cc}
D_{1}R_{1}^{-1}D_{1}^{\top }+2\lambda D_{2}R_{2}^{-1}D_{2}^{\top } &
-D_{2}R_{2}^{-1}D_{2}^{\top } \\
-D_{2}R_{2}^{-1}D_{2}^{\top } & 2\mu D_{2}R_{2}^{-1}D_{2}^{\top }%
\end{array}%
\right) , \\
\bar{Q}_{1} &=&\frac{1}{1+4\lambda \mu }\left(
\begin{array}{cc}
Q_{1}+2\mu Q_{2} & -Q_{2} \\
-Q_{2} & 2\lambda Q_{2}%
\end{array}%
\right) , \\
\bar{\Phi}_{1} &=&\frac{1}{1+4\lambda \mu }\left(
\begin{array}{cc}
\Phi _{1}+2\mu \Phi _{2} & -\Phi _{2} \\
-\Phi _{2} & 2\lambda \Phi _{2}%
\end{array}%
\right) .
\end{eqnarray*}%
Now it is easy to check that $\mathcal{\bar{B}}_{1},\mathcal{\bar{B}}_{2},%
\mathcal{\bar{D}}_{1},\mathcal{\bar{D}}_{2},\bar{Q}_{1}$ and $\bar{\Phi}_{1}$
are symmetric and positive definite. Repeating the approach above, we can
derive a standard backward stochastic Riccati equation as follows:
\begin{eqnarray}
\mathrm{d}\mathcal{\bar{R}}\left( t\right) &=&[\mathcal{\bar{A}}\left(
t\right) \mathcal{\bar{R}}\left( t\right) +\mathcal{\bar{C}}^{\top }\left(
t\right) \bar{\Xi}\left( t\right) +\bar{Q}_{1}\left( t\right)  \notag \\
&&+\mathcal{\bar{R}}\left( t\right) \left[ \mathcal{\bar{A}}\left( t\right) -%
\mathcal{\bar{B}}_{1}\left( t\right) \mathcal{\bar{R}}\left( t\right) -%
\mathcal{\bar{B}}_{2}\left( t\right) \bar{\Xi}\left( t\right) \right]  \notag
\\
&&+\bar{\Psi}\left( t\right) \left[ \mathcal{\bar{C}}\left( t\right) -%
\mathcal{\bar{D}}_{1}\left( t\right) \mathcal{\bar{R}}\left( t\right) -%
\mathcal{\bar{D}}_{2}\left( t\right) \bar{\Xi}\left( t\right) \right] ]%
\mathrm{d}t  \notag \\
&&-\bar{\Psi}\left( t\right) \mathrm{d}W\left( t\right)  \notag \\
\mathcal{\bar{R}}\left( T\right) &=&\bar{\Phi}_{1}.  \label{sriccati}
\end{eqnarray}%
where
\begin{equation*}
\bar{\Xi}\left( t\right) =\left( I+\mathcal{\bar{R}}\left( t\right) \mathcal{%
\bar{D}}_{2}\left( t\right) \right) ^{-1}\left[ \mathcal{\bar{R}}\left(
t\right) \mathcal{\bar{C}}\left( t\right) -\mathcal{\bar{R}}\left( t\right)
\mathcal{\bar{D}}_{1}\left( t\right) \mathcal{\bar{R}}\left( t\right) +\bar{%
\Psi}\left( t\right) \right] .
\end{equation*}%
From Tang \cite{Tang03}, we know that the Riccati equation (\ref{sriccati})
admits a unique solution. Moreover, we have
\begin{eqnarray*}
\bar{P} &=&\mathcal{\bar{R}}\bar{X}, \\
\bar{Q} &=&\bar{\Xi}\bar{X}, \\
\mathcal{R} &=&\Upsilon \mathcal{\bar{R}}\text{,} \\
\Psi &=&\Upsilon \bar{\Psi}.
\end{eqnarray*}%
From the fact
\begin{eqnarray*}
P &=&\Upsilon \bar{P}=\Upsilon \mathcal{\bar{R}}\bar{X}=\Upsilon \mathcal{%
\bar{R}}X, \\
Q &=&\Upsilon \bar{Q}=\Upsilon \bar{\Xi}\bar{X}=\Upsilon \bar{\Xi}X,
\end{eqnarray*}%
we state that the AOL Stackelberg solution $\left( u^{\ast }(\cdot ),v^{\ast
}(\cdot )\right) $ presents a feedback representation with respect to state $%
\left( x,k\right) $. The proof is completed. ~\hfill $\Box $

\subsection{The ACLM fashion}

In this subsection, we aforehand suppose that the derivative $\frac{\partial
u}{\partial x}$ to be bounded. As we shall see that the derivative enters
into the coefficient of the adjoint equation, its boundedness ensures the
well-posedness of the leader's problem whenever affine strategies are
considered. For simplicity, we study a one-dimensional linear quadratic
game, with the state equation as follows:%
\begin{equation}
\left\{
\begin{array}{rcl}
\mathrm{d}x & = & \left[ Ax+B_{1}u+B_{2}v\right] \mathrm{d}t+\left[
Cx+D_{1}u+D_{2}v\right] \mathrm{d}W\left( t\right) , \\
x\left( 0\right) & = & x_{0}\in \mathbb{R}.%
\end{array}%
\right.  \label{aclms}
\end{equation}%
The cost functionals for the leader and the follower to minimize are given,
respectively, as follows:%
\begin{equation}
\mathcal{J}_{1}\left( u,v\right) =\frac{1}{2}\mathbb{E}\left[ \int_{0}^{T}%
\left[ \left( Q_{1}x^{2}+R_{1}u^{2}\right) \right] \mathrm{d}t+\Phi
_{1}x^{2}\left( T\right) \right]  \label{aclml}
\end{equation}%
and
\begin{equation}
\mathcal{J}_{2}\left( u,v\right) =\frac{1}{2}\mathbb{E}\left[
\int_{0}^{T}\left( Q_{2}x^{2}+R_{2}v^{2}\right) \mathrm{d}t+\Phi
_{2}x^{2}\left( T\right) \right] .  \label{aclmf}
\end{equation}%
The admissible strategy spaces for the leader and the follower are denoted
by
\begin{eqnarray*}
\mathcal{U} &=&\Big \{u|u:\Omega \times \left[ 0,T\right] \times \mathbb{R}%
\rightarrow \mathbb{R}\text{ is }\mathcal{F}_{t}\text{-adapted for any }x\in
\mathbb{R},\text{ }u\left( t,x\right) \text{ is } \\
&&\text{continuously differentiable in }x\text{ for any }\left( t,x\right)
\in \Omega \times \left[ 0,T\right] \text{ satisfying }\left\vert \frac{%
\partial u}{\partial x}\right\vert <K\Big \}, \\
\mathcal{V} &=&\left\{ v\left\vert :\Omega \times \left[ 0,T\right] \times
\mathbb{R}\times \mathcal{U}\rightarrow \Gamma _{2}\text{ is }\mathcal{F}_{t}%
\text{-adapted for any }x\in \mathbb{R}\text{ and }u\in \mathcal{U}\right.
\right\} .
\end{eqnarray*}%
For any given $u\in \mathcal{U},$ the follower responses a unique optimal
strategy $v^{\ast }\in \mathcal{V}$. From (\ref{ACLM-f}), we get
\begin{equation*}
v^{\ast }\left( u\right) =\mathbf{P}_{\Gamma _{2}}[-R_{2}^{-1}(B_{2}^{\top
}p_{2}+D_{2}^{\top }q_{2})],
\end{equation*}%
where $p_{2}$ and $q_{2}$ satisfy
\begin{equation}
\left\{
\begin{array}{rcl}
-\mathrm{d}p_{2} & = & \left[ \left( A+B_{1}\frac{\partial u}{\partial x}%
\right) p_{2}+\left( C+D_{1}\frac{\partial u}{\partial x}\right) q_{2}+Q_{2}x%
\right] \mathrm{d}t-q_{2}\mathrm{d}W\left( t\right) , \\
p_{2}\left( T\right) & = & \Phi _{2}x\left( T\right) .%
\end{array}%
\right.  \label{adjaclm}
\end{equation}%
Now we formulate the leader's optimal control problem:%
\begin{equation*}
\min_{u\in \mathcal{U}}\mathcal{J}_{1}\left( u,v^{\ast }\right)
\end{equation*}%
subject to
\begin{equation}
\left\{
\begin{array}{rcl}
\mathrm{d}x & = & \left[ Ax+B_{1}u-B_{2}\mathbf{P}_{\Gamma
_{2}}[R_{2}^{-1}(B_{2}^{\top }p_{2}+D_{2}^{\top }q_{2})]\right] \mathrm{d}t
\\
&  & +\left[ Cx+D_{1}u-D_{2}\mathbf{P}_{\Gamma _{2}}[R_{2}^{-1}(B_{2}^{\top
}p_{2}+D_{2}^{\top }q_{2})]\right] \mathrm{d}W\left( t\right) , \\
-\mathrm{d}p_{2} & = & \left[ \left( A+B_{1}\frac{\partial u}{\partial x}%
\right) p_{2}+\left( C+D_{1}\frac{\partial u}{\partial x}\right) q_{2}+Q_{2}x%
\right] \mathrm{d}t-q_{2}\mathrm{d}W\left( t\right) , \\
p_{2}\left( T\right) & = & \Phi _{2}x\left( T\right) .%
\end{array}%
\right.  \label{aclmleader}
\end{equation}%
Note that FBSDEs (\ref{aclmleader}) are fully coupled. Due to the
boundedness of $\frac{\partial u}{\partial x},$ (\ref{aclmleader}) admit a
unique solution. Hence the leader's problem is well-posed. By means of
Theorem \ref{the2}, the leader can select his strategy among affine functions%
\begin{equation*}
u\left( t,x\right) =u_{2}\left( t\right) x+u_{1}\left( t\right) ,
\end{equation*}%
where $u_{2}$ and $u_{1}$ are adapted processes with $\left\vert
u_{2}\right\vert \leq K.$ Thus, the leader's problem can be described as
follows: The state equation is
\begin{equation}
\left\{
\begin{array}{rcl}
\mathrm{d}x & = & \left[ \left( A+B_{1}u_{2}\right) x+B_{1}u_{1}-B_{2}%
\mathbf{P}_{\Gamma _{2}}[R_{2}^{-1}(B_{2}^{\top }p_{2}+D_{2}^{\top }q_{2})]%
\right] \mathrm{d}t \\
&  & +\left[ \left( C+D_{1}u_{2}\right) x+D_{1}u_{1}-D_{2}\mathbf{P}_{\Gamma
_{2}}[R_{2}^{-1}(B_{2}^{\top }p_{2}+D_{2}^{\top }q_{2})]\right] \mathrm{d}%
W\left( t\right) , \\
-\mathrm{d}p_{2} & = & \left[ \left( A+B_{1}u_{2}\right) p_{2}+\left(
C+D_{1}u_{2}\right) q_{2}+Q_{2}x\right] \mathrm{d}t-q_{2}\mathrm{d}W\left(
t\right) , \\
p_{2}\left( T\right) & = & \Phi _{2}x\left( T\right) .%
\end{array}%
\right.  \label{aclmadj2}
\end{equation}%
The cost functional to be minimized:%
\begin{equation}
\min_{u\in \mathcal{U}}\mathcal{J}_{1}\left( u,v\right) =\min_{u\in \mathcal{%
U}}\frac{1}{2}\mathbb{E}\left[ \int_{0}^{T}\left[ \left(
Q_{1}x^{2}+R_{1}\left( u_{2}x+u_{1}\right) ^{2}\right) \right] \mathrm{d}%
t+\Phi _{1}x^{2}\left( T\right) \right] .  \label{aclmcost}
\end{equation}%
For any $\left( u_{1},u_{2}\right) $, we can get the existence and
uniqueness of the solution of (\ref{aclmadj2}) by the monotonicity condition
(Proposition \ref{mono}). Therefore, the leader's problem with strategies
restricted to be of affine form is well-posed. Now we apply the maximum
principle to obtain the Hamiltonian system and the related Riccati equation
for the leader's problem (\ref{aclmadj2})-(\ref{aclmcost}).

Set
\begin{eqnarray}
&&\mathcal{H}_{3}\left( t,u_{1},u_{2},x,\chi ,p_{1},q_{1},p_{2},q_{2}\right)
\notag \\
&=&p_{1}\left[ \left( A+B_{1}u_{2}\right) x+B_{1}u_{1}-B_{2}\mathbf{P}%
_{\Gamma _{2}}[R_{2}^{-1}(B_{2}^{\top }p_{2}+D_{2}^{\top }q_{2})]\right]
\notag \\
&&+q_{1}\left[ \left( C+D_{1}u_{2}\right) x+D_{1}u_{1}-D_{2}\mathbf{P}%
_{\Gamma _{2}}[R_{2}^{-1}(B_{2}^{\top }p_{2}+D_{2}^{\top }q_{2})]\right]
\notag \\
&&-\chi \left[ \left( A+B_{1}u_{2}\right) p_{2}+\left( C+D_{1}u_{2}\right)
q_{2}+Q_{2}x\right]  \notag \\
&&+\frac{1}{2}Q_{1}x^{2}+\frac{1}{2}R_{1}\left( u_{2}x+u_{1}\right) ^{2}.
\label{aclmh}
\end{eqnarray}%
Clearly, $\mathcal{H}_{3}$ is quadratic with respect to $u_{1}.$ Thus fix $%
u_{2}$ first, the minimizer of $\mathcal{H}_{3}$ can be obtained%
\begin{equation}
u_{1}^{\ast }=-u_{2}x-R_{1}^{-1}p_{1}B_{1}-R_{1}^{-1}q_{1}D_{1}.
\label{aclmof}
\end{equation}%
Inserting (\ref{aclmof}) into (\ref{aclmh}), we can see that the only term
involving $u_{2}$ is $-\chi \left( B_{1}p_{2}+D_{1}q_{2}\right) u_{2}.$
Hence, the optimal $u_{2}$ can be attained by
\begin{equation}
u_{2}^{\ast }=sgn\left( \chi \left( B_{1}p_{2}+D_{1}q_{2}\right) \right) K,
\label{aclmlcon}
\end{equation}%
where
\begin{equation*}
sgn\left( x\right) =\left\{
\begin{array}{cc}
1 & \text{if }x>0, \\
0 & \text{if }x=0, \\
-1 & \text{if }x<0.%
\end{array}%
\right.
\end{equation*}%
Now for simplicity, let $\Gamma _{2}=\mathbb{R}$. In this case,
\begin{eqnarray}
&&\mathcal{H}_{3}\left( t,u_{1},u_{2},x,\chi ,p_{1},q_{1},p_{2},q_{2}\right)
\notag \\
&=&p_{1}\left[ \left( A+B_{1}u_{2}\right)
x+B_{1}u_{1}-B_{2}R_{2}^{-1}(B_{2}^{\top }p_{2}+D_{2}^{\top }q_{2})]\right]
\notag \\
&&+q_{1}\left[ \left( C+D_{1}u_{2}\right)
x+D_{1}u_{1}-D_{2}R_{2}^{-1}(B_{2}^{\top }p_{2}+D_{2}^{\top }q_{2})\right]
\notag \\
&&-\chi \left[ \left( A+B_{1}u_{2}\right) p_{2}+\left( C+D_{1}u_{2}\right)
q_{2}+Q_{2}x\right]  \notag \\
&&+\frac{1}{2}Q_{1}x^{2}+\frac{1}{2}R_{1}\left( u_{2}x+u_{1}\right) ^{2}.
\end{eqnarray}%
So
\begin{equation*}
\left\{
\begin{array}{rcl}
\mathrm{d}x & = & \left[ \left( A+B_{1}u_{2}^{\ast }\right)
x+B_{1}u_{1}^{\ast }-B_{2}R_{2}^{-1}(B_{2}p_{2}+D_{2}^{\top }q_{2})\right]
\mathrm{d}t \\
&  & +\left[ \left( C+D_{1}u_{2}^{\ast }\right) x+D_{1}u_{1}^{\ast
}-D_{2}R_{2}^{-1}(B_{2}^{\top }p_{2}+D_{2}^{\top }q_{2})\right] \mathrm{d}%
W\left( t\right) , \\
\mathrm{d}\chi & = & \left[ \left( A+B_{1}u_{2}^{\ast }\right) \chi
+B_{2}R_{2}^{-1}B_{2}^{\top }p_{1}+D_{2}R_{2}^{-1}B_{2}^{\top }q_{1}\right]
\mathrm{d}t \\
&  & +\left[ \left( C+D_{1}u_{2}^{\ast }\right) \chi
+B_{2}R_{2}^{-1}D_{2}^{\top }p_{1}+D_{2}R_{2}^{-1}D_{2}^{\top }q_{1}\right]
\mathrm{d}W\left( t\right) , \\
-\mathrm{d}p_{1} & = & \left[ \left( A+B_{1}u_{2}^{\ast }\right)
p_{1}+\left( C+D_{1}u_{2}^{\ast }\right) q_{1}-\chi Q_{2}+Q_{1}x+R_{1}\left(
u_{2}^{\ast }x+u_{1}^{\ast }\right) u_{2}^{\ast }\right] \mathrm{d}t \\
&  & -q_{1}\left( t\right) \mathrm{d}W\left( t\right) , \\
x\left( 0\right) & = & x_{0},\text{ }\chi \left( 0\right) =0,p_{1}\left(
T\right) =-\Phi _{2}\chi \left( T\right) +\Phi _{1}x\left( T\right) .%
\end{array}%
\right.
\end{equation*}%
Put (\ref{aclmof}) and (\ref{aclmlcon}) into above, we get%
\begin{equation}
\left\{
\begin{array}{rcl}
\mathrm{d}x & = & \left[ Ax-B_{1}^{2}R_{1}^{-1}\left(
B_{1}p_{1}+D_{1}q_{1}\right) -B_{2}R_{2}^{-1}(B_{2}p_{2}+D_{2}^{\top }q_{2})%
\right] \mathrm{d}t \\
&  & +\left[ Cx-D_{1}^{2}R_{1}^{-1}\left( B_{1}p_{1}+D_{1}q_{1}\right)
-D_{2}R_{2}^{-1}(B_{2}p_{2}+D_{2}^{\top }q_{2})\right] \mathrm{d}W\left(
t\right) , \\
\mathrm{d}\chi & = & \left[ A\chi +B_{1}sgn\left( \chi \left(
B_{1}p_{2}+D_{1}q_{2}\right) \right) K\chi +B_{2}R_{2}^{-1}B_{2}^{\top
}p_{1}+D_{2}R_{2}^{-1}B_{2}^{\top }q_{1}\right] \mathrm{d}t \\
&  & +\left[ C\chi +D_{1}sgn\left( \chi \left( B_{1}p_{2}+D_{1}q_{2}\right)
\right) K\chi +B_{2}R_{2}^{-1}D_{2}^{\top }p_{1}+D_{2}R_{2}^{-1}D_{2}^{\top
}q_{1}\right] \mathrm{d}W\left( t\right) , \\
-\mathrm{d}p_{1} & = & \left[ Ap_{1}+Cq_{1}-\chi Q_{2}+Q_{1}x\right] \mathrm{%
d}t-q_{1}\left( t\right) \mathrm{d}W\left( t\right) , \\
-\mathrm{d}p_{2} & = & \big [\left( A+B_{1}sgn\left( \chi \left(
B_{1}p_{2}+D_{1}q_{2}\right) \right) K\right) p_{2} \\
&  & +\left( C+D_{1}sgn\left( \chi \left( B_{1}p_{2}+D_{1}q_{2}\right)
\right) K\right) q_{2}+Q_{2}x\big ]\mathrm{d}t-q_{2}\mathrm{d}W\left(
t\right) , \\
x\left( 0\right) & = & x_{0},\text{ }\chi \left( 0\right) =0,p_{1}\left(
T\right) =-\Phi _{2}\chi \left( T\right) +\Phi _{1}x\left( T\right)
,p_{2}\left( T\right) =\Phi _{2}x\left( T\right) .%
\end{array}%
\right.  \label{aclmhal}
\end{equation}%
Clearly, the Hamilton system (\ref{aclmhal}) is highly complex due to the
nonlinear term $sgn\left( \chi \left( B_{1}p_{2}+D_{1}q_{2}\right) \right) .$
Specifically, the state $\left( x,\chi \right) ^{\top }$ cannot be written
as a linear equation with respect to $\left( p_{1},p_{2}\right) ^{\top }.$
However, $sgn\left( \chi \left( B_{1}p_{2}+D_{1}q_{2}\right) \right) $
doesn't appear into $x,$ besides, the terminal condition $p_{1}\left(
T\right) ,$ $p_{2}\left( T\right) $ contain $x\left( T\right) .$ Hence, we
only regard $x$ as the state and postulate%
\begin{eqnarray}
\chi \left( t\right) &=&\alpha \left( t\right) x\left( t\right) ,  \label{r1}
\\
p_{1}\left( t\right) &=&\beta \left( t\right) x\left( t\right) ,  \label{r2}
\\
p_{2}\left( t\right) &=&\gamma \left( t\right) x\left( t\right) ,  \label{r3}
\end{eqnarray}%
where
\begin{eqnarray*}
\mathrm{d}\alpha \left( t\right) &=&\alpha _{1}\left( t\right) \mathrm{d}%
t+\alpha _{2}\left( t\right) \mathrm{d}W\left( t\right) , \\
\mathrm{d}\beta \left( t\right) &=&\beta _{1}\left( t\right) \mathrm{d}%
t+\beta _{2}\left( t\right) \mathrm{d}W\left( t\right) , \\
\mathrm{d}\gamma \left( t\right) &=&\gamma _{1}\left( t\right) \mathrm{d}%
t+\gamma _{2}\left( t\right) \mathrm{d}W\left( t\right) .
\end{eqnarray*}%
We first look for $q_{1},q_{2}$. Applying It\^{o}'s formula to $\beta \left(
t\right) x\left( t\right) $ and $\gamma \left( t\right) x\left( t\right) ,$
we have%
\begin{eqnarray*}
&&\beta \left[ Ax-B_{1}^{2}R_{1}^{-1}\left( B_{1}\beta x+D_{1}q_{1}\right)
-B_{2}R_{2}^{-1}(B_{2}\gamma x+D_{2}^{\top }q_{2})\right] \mathrm{d}t \\
&&+\beta \left[ Cx-D_{1}^{2}R_{1}^{-1}\left( B_{1}\beta x+D_{1}q_{1}\right)
-D_{2}R_{2}^{-1}(B_{2}\gamma x+D_{2}^{\top }q_{2})\right] \mathrm{d}W\left(
t\right) \\
&&+x\beta _{1}\mathrm{d}t+x\beta _{2}\mathrm{d}W\left( t\right) \\
&&+\beta _{2}\left[ Cx-D_{1}^{2}R_{1}^{-1}\left( B_{1}\beta
x+D_{1}q_{1}\right) -D_{2}R_{2}^{-1}(B_{2}\gamma x+D_{2}^{\top }q_{2})\right]
\mathrm{d}t \\
&=&\mathrm{d}p_{1} \\
&=&-\left[ A\beta x+Cq_{1}-\alpha xQ_{2}+Q_{1}x\right] \mathrm{d}t+q_{1}%
\mathrm{d}W\left( t\right) ,
\end{eqnarray*}%
and
\begin{eqnarray*}
&&\gamma \left[ Ax-B_{1}^{2}R_{1}^{-1}\left( B_{1}\beta x+D_{1}q_{1}\right)
-B_{2}R_{2}^{-1}(B_{2}\gamma x+D_{2}^{\top }q_{2})\right] \mathrm{d}t \\
&&+\gamma \left[ Cx-D_{1}^{2}R_{1}^{-1}\left( B_{1}\beta x+D_{1}q_{1}\right)
-D_{2}R_{2}^{-1}(B_{2}\gamma x+D_{2}^{\top }q_{2})\right] \mathrm{d}W\left(
t\right) \\
&&+x\gamma _{1}\mathrm{d}t+x\gamma _{2}\mathrm{d}W\left( t\right) \\
&=&\mathrm{d}p_{2} \\
&=&-\big [\left( A+B_{1}sgn\left( \alpha x\left(
B_{1}p_{2}+D_{1}q_{2}\right) \right) K\right) p_{2} \\
&&+\left( C+D_{1}sgn\left( \alpha x\left( B_{1}p_{2}+D_{1}q_{2}\right)
\right) K\right) q_{2}+Q_{2}x\big ]\mathrm{d}t+q_{2}\mathrm{d}W\left(
t\right) ,
\end{eqnarray*}%
So%
\begin{eqnarray*}
q_{1} &=&\beta Cx-\beta D_{1}^{2}R_{1}^{-1}B_{1}\beta x-\beta
D_{1}^{2}R_{1}^{-1}D_{1}q_{1} \\
&&-\beta D_{2}R_{2}^{-1}B_{2}\gamma x-\beta D_{2}R_{2}^{-1}D_{2}^{\top
}q_{2}+x\beta _{2}, \\
q_{2} &=&\gamma Cx-\gamma D_{1}^{2}R_{1}^{-1}B_{1}\beta x-\gamma
D_{1}^{2}R_{1}^{-1}D_{1}q_{1} \\
&&-\gamma D_{2}R_{2}^{-1}B_{2}\gamma x-\gamma D_{2}R_{2}^{-1}D_{2}^{\top
}q_{2}+x\gamma _{2},
\end{eqnarray*}%
namely,%
\begin{equation*}
\left\{
\begin{array}{l}
\left( 1+\beta D_{1}^{2}R_{1}^{-1}D_{1}\right) q_{1}+\beta
D_{2}R_{2}^{-1}D_{2}^{\top }q_{2}=\Xi _{1}x, \\
\gamma D_{1}^{2}R_{1}^{-1}D_{1}q_{1}+\left( 1+\gamma
D_{2}R_{2}^{-1}D_{2}^{\top }\right) q_{2}=\Xi _{2}x,%
\end{array}%
\right.
\end{equation*}%
where
\begin{eqnarray*}
\Xi _{1} &=&\beta C-\beta D_{1}^{2}R_{1}^{-1}B_{1}\beta -\beta
D_{2}R_{2}^{-1}B_{2}\gamma +\beta _{2}, \\
\Xi _{2} &=&\gamma C-\gamma D_{1}^{2}R_{1}^{-1}B_{1}\beta -\gamma
D_{2}R_{2}^{-1}B_{2}\gamma +\gamma _{2}.
\end{eqnarray*}%
It is easy to get%
\begin{eqnarray*}
q_{1} &=&\Delta _{1}x, \\
q_{2} &=&\Delta _{2}x,
\end{eqnarray*}%
where
\begin{eqnarray*}
\Delta _{1} &=&\frac{\left[ \Xi _{1}\left( 1+\gamma
D_{2}R_{2}^{-1}D_{2}^{\top }\right) -\Xi _{2}\beta
D_{2}R_{2}^{-1}D_{2}^{\top }\right] }{\left( 1+\beta
D_{1}^{2}R_{1}^{-1}D_{1}\right) \left( 1+\gamma D_{2}R_{2}^{-1}D_{2}^{\top
}\right) -\beta D_{2}R_{2}^{-1}D_{2}^{\top }\gamma D_{1}^{2}R_{1}^{-1}D_{1}},
\\
\Delta _{2} &=&\frac{\left[ \Xi _{2}\left( 1+\beta
D_{1}^{2}R_{1}^{-1}D_{1}\right) -\Xi _{1}\gamma D_{1}^{2}R_{1}^{-1}D_{1}%
\right] }{\left( 1+\beta D_{1}^{2}R_{1}^{-1}D_{1}\right) \left( 1+\gamma
D_{2}R_{2}^{-1}D_{2}^{\top }\right) -\beta D_{2}R_{2}^{-1}D_{2}^{\top
}\gamma D_{1}^{2}R_{1}^{-1}D_{1}}.
\end{eqnarray*}%
Moreover,
\begin{eqnarray*}
\beta _{1} &=&-\left[ A\beta +C\Delta _{1}-\alpha Q_{2}+Q_{1}\right] \\
&&-\beta \left[ A-B_{1}^{2}R_{1}^{-1}\left( B_{1}\beta +D_{1}\Delta
_{1}\right) -B_{2}R_{2}^{-1}(B_{2}\gamma +D_{2}^{\top }\Delta _{2})\right] \\
&&-\beta _{2}\left[ C-D_{1}^{2}R_{1}^{-1}\left( B_{1}\beta +D_{1}\Delta
_{1}\right) -D_{2}R_{2}^{-1}(B_{2}\gamma +D_{2}^{\top }\Delta _{2})\right]
\end{eqnarray*}%
and%
\begin{eqnarray*}
\gamma _{1} &=&-\big [\left( A+B_{1}sgn\left( \alpha \left( B_{1}\gamma
+D_{1}\Delta _{2}\right) \right) K\right) \gamma \\
&&+\left( C+D_{1}sgn\left( \alpha \left( B_{1}\gamma +D_{1}\Delta
_{2}\right) \right) K\right) \Delta _{2}+Q_{2}\big ] \\
&&-\gamma \left[ A-B_{1}^{2}R_{1}^{-1}\left( B_{1}\beta +D_{1}\Delta
_{1}\right) -B_{2}R_{2}^{-1}(B_{2}\gamma +D_{2}^{\top }\Delta _{2})\right]
\mathrm{.}
\end{eqnarray*}%
Repeating the method used above, we have%
\begin{eqnarray*}
&&\alpha \left[ Ax-B_{1}^{2}R_{1}^{-1}\left( B_{1}\beta x+D_{1}q_{1}\right)
-B_{2}R_{2}^{-1}(B_{2}\gamma x+D_{2}^{\top }\rho _{2}x)\right] \mathrm{d}%
t+x\alpha _{1}\mathrm{d}t+x\alpha _{2}\mathrm{d}W \\
&&+\alpha \left[ Cx-D_{1}^{2}R_{1}^{-1}\left( B_{1}\beta x+D_{1}q_{1}\right)
-D_{2}R_{2}^{-1}(B_{2}\gamma x+D_{2}^{\top }q_{2})\right] \mathrm{d}W\left(
t\right) \\
&&+\alpha _{2}\left[ Cx-D_{1}^{2}R_{1}^{-1}\left( B_{1}\beta
x+D_{1}q_{1}\right) -D_{2}R_{2}^{-1}(B_{2}\gamma x+D_{2}^{\top }q_{2})\right]
\mathrm{d}t \\
&=&\mathrm{d}\chi \\
&=&\left[ A\alpha x+B_{1}sgn\left( \alpha x\left( B_{1}\gamma
x+D_{1}q_{2}\right) \right) K\alpha x+B_{2}R_{2}^{-1}B_{2}^{\top }\beta
x+D_{2}R_{2}^{-1}B_{2}^{\top }q_{1}\right] \mathrm{d}t \\
&&+\left[ C\alpha x+D_{1}sgn\left( \alpha x\left( B_{1}\gamma
x+D_{1}q_{2}\right) \right) K\alpha x+B_{2}R_{2}^{-1}D_{2}^{\top }\beta
x+D_{2}R_{2}^{-1}D_{2}^{\top }q_{1}\right] \mathrm{d}W\left( t\right) .
\end{eqnarray*}%
By comparing (\ref{aclmhal}) and (\ref{r1}), we have%
\begin{eqnarray*}
\alpha _{2} &=&\alpha \left[ D_{1}^{2}R_{1}^{-1}\left( B_{1}\beta
+D_{1}\Delta _{1}\right) +D_{2}R_{2}^{-1}(B_{2}\gamma +D_{2}^{\top }\Delta
_{2})\right] \\
&&+D_{1}sgn\left( \alpha \left( B_{1}\gamma +D_{1}\Delta _{2}\right) \right)
\alpha +B_{2}R_{2}^{-1}D_{2}^{\top }\beta +D_{2}R_{2}^{-1}D_{2}^{\top
}\Delta _{1}, \\
\alpha _{1} &=&B_{1}sgn\left( \alpha \left( B_{1}\gamma +D_{1}\Delta
_{2}\right) \right) K\alpha +B_{2}R_{2}^{-1}B_{2}^{\top }\beta
+D_{2}R_{2}^{-1}B_{2}^{\top }\Delta _{1} \\
&&+\alpha \left[ B_{1}^{2}R_{1}^{-1}\left( B_{1}\beta +D_{1}\Delta
_{1}\right) +B_{2}R_{2}^{-1}(B_{2}\gamma +D_{2}^{\top }\Delta _{2})\right] \\
&&+\alpha _{2}\left[ C+D_{1}^{2}R_{1}^{-1}\left( B_{1}\beta +D_{1}\Delta
_{1}\right) +D_{2}R_{2}^{-1}(B_{2}\gamma +D_{2}^{\top }\Delta _{2})\right] .
\end{eqnarray*}%
Now we are able to announce a new system as follows:%
\begin{equation}
\left\{
\begin{array}{rcl}
\mathrm{d}x^{\ast } & = & \left[ A-B_{1}^{2}R_{1}^{-1}\left( B_{1}\beta
+D_{1}\Delta _{1}\right) -B_{2}R_{2}^{-1}(B_{2}\gamma +D_{2}^{\top }\Delta
_{2})\right] x^{\ast }\mathrm{d}t \\
&  & +\left[ C-D_{1}^{2}R_{1}^{-1}\left( B_{1}\beta +D_{1}\Delta _{1}\right)
-D_{2}R_{2}^{-1}(B_{2}\gamma +D_{2}^{\top }\Delta _{2})\right] x^{\ast }%
\mathrm{d}W\left( t\right) , \\
\mathrm{d}\alpha & = & \big [B_{1}sgn\left( \alpha \left( B_{1}\gamma
+D_{1}\Delta _{2}\right) \right) K\alpha +B_{2}R_{2}^{-1}B_{2}^{\top }\beta
+D_{2}R_{2}^{-1}B_{2}^{\top }\Delta _{1} \\
&  & +\alpha \left[ B_{1}^{2}R_{1}^{-1}\left( B_{1}\beta +D_{1}\Delta
_{1}\right) +B_{2}R_{2}^{-1}(B_{2}\gamma +D_{2}^{\top }\Delta _{2})\right]
\\
&  & +\alpha _{2}\left[ C+D_{1}^{2}R_{1}^{-1}\left( B_{1}\beta +D_{1}\Delta
_{1}\right) +D_{2}R_{2}^{-1}(B_{2}\gamma +D_{2}^{\top }\Delta _{2})\right] %
\big ]\mathrm{d}t \\
&  & \big [\alpha \left[ D_{1}^{2}R_{1}^{-1}\left( B_{1}\beta +D_{1}\Delta
_{1}\right) +D_{2}R_{2}^{-1}(B_{2}\gamma +D_{2}^{\top }\Delta _{2})\right]
\\
&  & +D_{1}sgn\left( \alpha \left( B_{1}\gamma +D_{1}\Delta _{2}\right)
\right) K\alpha +B_{2}R_{2}^{-1}D_{2}^{\top }\beta
+D_{2}R_{2}^{-1}D_{2}^{\top }\Delta _{1}\big ]\mathrm{d}W\left( t\right) ,
\\
\mathrm{d}\beta & = & \big \{-\left[ A\beta +C\Delta _{1}-\alpha Q_{2}+Q_{1}%
\right] \\
&  & -\beta \left[ A-B_{1}^{2}R_{1}^{-1}\left( B_{1}\beta +D_{1}\Delta
_{1}\right) -B_{2}R_{2}^{-1}(B_{2}\gamma +D_{2}^{\top }\Delta _{2})\right]
\\
&  & -\beta _{2}\left[ C-D_{1}^{2}R_{1}^{-1}\left( B_{1}\beta +D_{1}\Delta
_{1}\right) -D_{2}R_{2}^{-1}(B_{2}\gamma +D_{2}^{\top }\Delta _{2})\right] %
\big \}\mathrm{d}t+\beta _{2}\mathrm{d}W\left( t\right) , \\
\mathrm{d}\gamma & = & \big \{-\big [\left( A+B_{1}sgn\left( \alpha \left(
B_{1}\gamma +D_{1}\Delta _{2}\right) \right) K\right) \gamma \\
&  & +\left( C+D_{1}sgn\left( \alpha \left( B_{1}\gamma +D_{1}\Delta
_{2}\right) \right) K\right) \Delta _{2}+Q_{2}\big ] \\
&  & -\gamma \left[ A-B_{1}^{2}R_{1}^{-1}\left( B_{1}\beta +D_{1}\Delta
_{1}\right) -B_{2}R_{2}^{-1}(B_{2}\gamma +D_{2}^{\top }\Delta _{2})\right] %
\big \}\mathrm{d}t+\gamma _{2}\mathrm{d}W\left( t\right) , \\
\alpha \left( 0\right) & = & 0,\text{ }\beta \left( T\right) =-\Phi
_{2}\alpha \left( T\right) +\Phi _{1},\gamma \left( T\right) =\Phi _{2}.%
\end{array}%
\right.  \label{aclmlarge}
\end{equation}%
Suppose that FBSDEs (\ref{aclmlarge}) admit a unique solution, denoted by $%
\left( x^{\ast },\alpha ,\beta ,\gamma ,\beta _{2},\gamma _{2}\right) ,$
which actually solve the Hamiltonian system (\ref{aclmhal}). As a result, a
candidate for the leader's optimal strategy can be expressed as
\begin{eqnarray}
u\left( t,x\right) &=&sgn\left( \alpha \left( B_{1}\gamma +D_{1}\Delta
_{2}\right) \right) Kx-sgn\left( \alpha \left( B_{1}\gamma +D_{1}\Delta
_{2}\right) \right) Kx^{\ast }\left( t\right)  \notag \\
&&-R_{1}^{-1}\beta B_{1}x^{\ast }\left( t\right) -R_{1}^{-1}\Delta
_{1}D_{1}x^{\ast }\left( t\right) .  \label{opleader}
\end{eqnarray}

\begin{remark}
From (\ref{opleader}), it indicates that, comparing with Bensoussan et al.
\cite{BCS2015}, the matrices $D_{1},D_{2}$ will impose on $u\left(
t,x\right) $ via $D_{1}\Delta _{2}$ and $\Delta _{1}D_{1}.$ If $%
D_{1}=D_{2}=0,$ (\ref{opleader}) reduces to Bensoussan et al.'s type (see
\cite{BCS2015}).
\end{remark}

\section{Conclusions and remarks}

\label{sec4}

In this paper, we are concerned on the solutions of stochastic Stackelberg
differential games within two information structures: AOL and ACLM patterns
under convex control domain. Having maximum principle for the former kind of
game as a basis, we give the necessary conditions for the leader's optimal
strategy in the latter game. To illustrate our theoretic results, we study
the LQ stochastic Stackelberg differential games. For AOL case, on the one
hand, we prove the existence and uniqueness of Hamiltonian system of
leader's with projection operator and derive a kind of standard backward
stochastic Riccati equation. For ACLM case, we also give a Riccati equation
with non-linear term and then provide the leader's optimal strategy.

There are some topics deserved to study displaying in the following: (i) The
Riccati equation derived in the framework of ACLM case consists of complex
coupled FBSDEs with non-linear term. The general conditions to guarantee the
existence and uniqueness is not known. (ii) As observed that the control set
is limited to convex, a natural question arises: How to establish the
maximum principle for general case, namely, non-convex control domain?
Certainly, the second-order adjoint equation is employed, which makes the
system extremely complicated. The Hamiltonian system in this situation
actually involves six types of (forward or backward) stochastic differential
equations. (iii) It is necessary to point out that the Stackelberg game
considered in this paper is limited to in the complete information
background. In other word, both the leader and the follower can observe the
state process of the stochastic system directly, which, however, is
impossible in reality. As a matter of fact, both of them can only announce
partial information because of the market competition, information-delay,
private information and limitation policy by the government, \textit{etc}.
Therefore, it is necessary to study the Stackelberg game with partial and
asymmetric information (see \cite{SWX2016, SWX2017}).\ These possible
extensions to the Stackelberg stochastic differential game no doubt promise
to be interesting research directions. We shall response these challenging
topics in our future work.

\appendix
\label{APP}

\section{Properties of projection $\mathbf{P}_{\Gamma }$}

For the readers' convenience, let us recall the following properties of
projection $\mathbf{P}_{\Gamma}$ onto a closed convex set, see \cite{Brezis
2010}, Chapter 5.

\begin{proposition}
\label{projection theorem} For a nonempty closed convex set $\Gamma\subset%
\mathbb{R}^m$, for every $x\in\mathbb{R}^m$, there exists a unique $%
x^\ast\in\Gamma$, such that
\begin{equation*}
|x-x^\ast|=\min_{y\in\Gamma}|x-y|=:dist(x,\Gamma).
\end{equation*}
Moreover, $x^\ast$ is characterized by the property
\begin{equation}  \label{projection characterization}
x^\ast\in\Gamma, \quad \big<x^\ast-x, x^\ast-y\big>\leq0 \qquad \forall y
\in \Gamma.
\end{equation}
The above element $x^\ast$ is called the projection of $x$ onto $\Gamma$ and
is denoted by $\mathbf{P}_{\Gamma}[x]$.
\end{proposition}

From above theorem, it is easy to show that

\begin{proposition}
\label{mono}Let $\Gamma \subset \mathbb{R}^{m}$ be a nonempty closed convex
set, then we have
\begin{equation}
\big|\mathbf{P}_{\Gamma }[x]-\mathbf{P}_{\Gamma }[y]\big|^{2}\leq \big<%
\mathbf{P}_{\Gamma }[x]-\mathbf{P}_{\Gamma }[y],x-y\big>.
\label{projection inequality}
\end{equation}
\end{proposition}


\begin{proposition}
\label{Plip}Let $\Gamma \subset \mathbb{R}^{m}$ be a nonempty closed convex
set, then the projection $\mathbf{P}_{\Gamma }$ does not increase the
distance, i.e.
\begin{equation}
\big|\mathbf{P}_{\Gamma }[x]-\mathbf{P}_{\Gamma }[y]\big|\leq \big|x-y\big|.
\label{projection lipschitz}
\end{equation}
\end{proposition}


Now let us consider $\mathbb{R}^m$ and the projection $\mathbf{P}_{\Gamma}$
both with the norm $\|\cdot\|_{R_0}:=\langle R_0^{\frac{1}{2}}\cdot, R_0^{%
\frac{1}{2}}\cdot\rangle$, from (\ref{projection inequality}), we have

\begin{proposition}
\label{mono2}Let $\Gamma \subset \mathbb{R}^{m}$ be a nonempty closed convex
set, then
\begin{equation*}
\langle \langle \mathbf{P}_{\Gamma }[x]-\mathbf{P}_{\Gamma }[y],x-y\rangle
\rangle =\left\langle R^{\frac{1}{2}}\bigg(\mathbf{P}_{\Gamma }[x]-\mathbf{P}%
_{\Gamma }[y]\bigg),R^{\frac{1}{2}}(x-y)\right\rangle \geq 0.
\end{equation*}
\end{proposition}

\section{The proof of Theorem \protect\ref{the1}}

\paragraph{Proof.}

(\textbf{Uniqueness}) Suppose that there exists two solutions: $%
(x^{1},p_{2}^{1},q_{2}^{1})$, $(x^{2},p_{2}^{2},q_{2}^{2})$ and denote%
\begin{equation*}
\hat{x}=x^{1}-x^{2},\quad \hat{p}_{2}=p_{2}^{1}-p_{2}^{2},\quad \hat{q}%
_{2}=q_{2}^{1}-q_{2}^{2}.
\end{equation*}%
Then, we have
\begin{equation}
\left\{
\begin{array}{rcl}
\mathrm{d}\hat{x}\left( t\right) & = & \left[ A\left( t\right) \hat{x}\left(
t\right) +B_{2}\left( t\right) \hat{\varphi}_{2}(t,\hat{p}_{2}\left(
t\right) ,\hat{q}_{2}\left( t\right) )\right] \mathrm{d}t \\
&  & +\left[ C\left( t\right) \hat{x}\left( t\right) +D_{2}\left( t\right)
\hat{\varphi}_{2}(t,\hat{p}_{2}\left( t\right) ,\hat{q}_{2}\left( t\right) )%
\right] \mathrm{d}W\left( t\right) , \\
-\mathrm{d}\hat{p}_{2}\left( t\right) & = & \left[ A^{\top }\left( t\right)
\hat{p}_{2}\left( t\right) +C^{\top }\left( t\right) \hat{q}_{2}\left(
t\right) -Q_{2}\left( t\right) \hat{x}\left( t\right) \right] \mathrm{d}t-%
\hat{q}_{2}\left( t\right) W\left( t\right) , \\
\hat{x}\left( 0\right) & = & 0,\text{ }\hat{p}_{2}\left( T\right) =-\Phi _{2}%
\hat{x}\left( T\right) .%
\end{array}%
\right.  \label{e221}
\end{equation}%
with
\begin{eqnarray*}
\hat{\varphi}_{2}(t,\hat{p}_{2}\left( t\right) ,\hat{q}_{2}\left( t\right) )
&=&\varphi _{2}(t,p_{2}^{1}\left( t\right) ,q_{2}^{1}\left( t\right)
)-\varphi _{2}(t,p_{2}^{2}\left( t\right) ,q_{2}^{2}\left( t\right) ) \\
&=&\mathbf{P}_{\Gamma _{2}}[R_{2}^{-1}(t)(B_{2}^{\top }(t)p_{2}^{1}\left(
t\right) +D_{2}^{\top }(t)q_{2}^{1}\left( t\right) )] \\
&&-\mathbf{P}_{\Gamma _{2}}[R_{2}^{-1}(t)(B_{2}^{\top }(t)p_{2}^{2}\left(
t\right) +D_{2}^{\top }(t)q_{2}^{2}\left( t\right) )]
\end{eqnarray*}%
First, applying It\^{o}'s formula to $\big<\hat{p},\hat{x}\big>$ and taking
expectations on both sides (noting the monotonicity property of $\widehat{%
\varphi },$ see Proposition \ref{mono2}$),$ we have:%
\begin{eqnarray*}
0 &=&\mathbb{E}\left\langle \Phi _{2}\hat{x}\left( T\right) ,\hat{x}\left(
T\right) \right\rangle \\
&&+\mathbb{E}\bigg [\int_{0}^{T}\left( \left\langle \left( B_{2}^{\top
}\left( t\right) \hat{p}_{2}\left( t\right) +D_{2}^{\top }\left( t\right)
\hat{q}_{2}\left( t\right) \right) ,\hat{\varphi}_{2}(t,\hat{p}_{2}\left(
t\right) ,\hat{q}_{2}\left( t\right) )\right\rangle +\left\langle \hat{x}%
\left( t\right) ,Q_{2}\left( t\right) \hat{x}\left( t\right) \right\rangle
\right) \mathrm{d}t\bigg ] \\
&\geq &\mathbb{E}\left\langle \Phi _{2}\hat{x}\left( T\right) ,\hat{x}\left(
T\right) \right\rangle +\mathbb{E}\bigg [\int_{0}^{T}\left\langle \hat{x}%
\left( t\right) ,Q_{2}\left( t\right) \hat{x}\left( t\right) \right\rangle
\mathrm{d}t\bigg ]
\end{eqnarray*}%
Thus, $\Phi _{2}\hat{x}\left( T\right) =0$ and $Q_{2}\left( t\right) \hat{x}%
\left( t\right) =0$ which implies $\hat{p}_{2}\left( t\right) \equiv 0,$ $%
\hat{q}_{2}\left( t\right) \equiv 0.$ Next, we have $\hat{\varphi}_{2}(t,%
\hat{p}_{2}\left( t\right) ,\hat{q}_{2}\left( t\right) )\equiv 0$ which
further implies $\hat{x}\left( t\right) \equiv 0.$ Hence the uniqueness
follows.

\vskip 12pt

(\textbf{Existence}) Consider a family of parameterized FBSDEs as follows%
\footnote{%
For simplicity, the dependence of coefficients on time variable $t$ is
suppressed.}:%
\begin{equation*}
\left\{
\begin{array}{rcl}
\mathrm{d}x^{\alpha } & = & \left[ \alpha \mathbb{B}\left( x^{\alpha
},p_{2}^{\alpha },q_{2}^{\alpha }\right) +\psi \right] \mathrm{d}t+\left[
\alpha \mathbb{C}\left( x^{\alpha },p_{2}^{\alpha },q_{2}^{\alpha }\right)
+\phi \right] \mathrm{d}W\left( t\right) , \\
-\mathrm{d}p_{2}^{\alpha } & = & \left[ \alpha \mathbb{F}\left( x^{\alpha
},p_{2}^{\alpha },q_{2}^{\alpha }\right) +\zeta \right] \mathrm{d}%
t-q_{2}^{\alpha }\mathrm{d}W\left( t\right) \\
x^{\alpha }\left( 0\right) & = & x_{0},p_{2}^{\alpha }\left( T\right)
=-\alpha \Phi _{2}\hat{x}\left( T\right) +\eta ,%
\end{array}%
\right.
\end{equation*}%
with%
\begin{equation*}
\left\{
\begin{array}{l}
\mathbb{B}\left( t,x,p_{2},q_{2}\right) =Ax+B_{2}\varphi _{2}(t,p_{2},q_{2})
\\
\mathbb{C}\left( t,x,p_{2},q_{2}\right) =Cx+D_{2}\left( t\right) \varphi
_{2}(t,p_{2},q_{2}) \\
\mathbb{F}\left( t,x,p_{2},q_{2}\right) =A^{\top }p_{2}+C^{\top }q_{2}-Q_{2}x%
\end{array}%
\right.
\end{equation*}%
Here, $(\psi ,\phi ,\zeta )$ are given processes in $\mathcal{M}^{2}(0,T;%
\mathbb{R}^{n})\times \mathcal{M}^{2}(0,T;\mathbb{R}^{n})\times \mathcal{M}%
^{2}(0,T;\mathbb{R}^{n}),$ and $\eta $ is a $\mathbb{R}^{n}$-valued square
integrable random variable which is $\mathcal{F}_{T}$-measurable. When $%
\alpha =0,$ we have a decoupled FBSDEs whose solvability is trivial:
\begin{equation*}
\left\{
\begin{array}{rcl}
\mathrm{d}x & = & \psi \mathrm{d}t+\phi \mathrm{d}W\left( t\right) , \\
-\mathrm{d}p_{2} & = & \zeta \mathrm{d}t-q_{2}\mathrm{d}W\left( t\right) \\
x\left( 0\right) & = & x_{0},p_{2}\left( T\right) =\eta ,%
\end{array}%
\right.
\end{equation*}%
Denote
\begin{equation*}
\widetilde{\mathcal{M}}(0,T)=\mathcal{M}^{2}(0,T;\mathbb{R}^{n})\times
\mathcal{M}^{2}(0,T;\mathbb{R}^{n})\times \mathcal{M}^{2}(0,T;\mathbb{R}%
^{n}).
\end{equation*}%
Now we introduce a mapping ${I}_{\alpha _{0}}:(x,p_{2},q_{2})\in \widetilde{%
\mathcal{M}}(0,T)\longrightarrow (X,P_{2},Q_{2})\in \widetilde{\mathcal{M}}%
(0,T)$ via the following FBSDEs:%
\begin{equation*}
\left\{
\begin{array}{rcl}
\mathrm{d}X & = & \left[ \alpha _{0}\mathbb{B}\left( X,P_{2},Q_{2}\right)
+\delta \mathbb{B}\left( x,p_{2},q_{2}\right) +\psi \right] \mathrm{d}t \\
&  & +\left[ \alpha _{0}\mathbb{C}\left( X,P_{2},Q_{2}\right) +\delta
\mathbb{C}\left( x,p_{2},q_{2}\right) +\phi \right] \mathrm{d}W\left(
t\right) , \\
-\mathrm{d}P_{2} & = & \left[ \alpha _{0}\mathbb{F}\left(
X,P_{2},Q_{2}\right) +\delta \mathbb{F}\left( x,p_{2},q_{2}\right) +\zeta %
\right] \mathrm{d}t-Q\mathrm{d}W\left( t\right) \\
X\left( 0\right) & = & x_{0},\text{ }P\left( T\right) =-\alpha _{0}\Phi
_{2}X\left( T\right) -\delta \Phi _{2}x\left( T\right) +\eta ,%
\end{array}%
\right.
\end{equation*}%
Considering ${I}_{\alpha _{0}}:(x,p_{2},q_{2})\longrightarrow
(X,P_{2},Q_{2}) $ and ${I}_{\alpha _{0}}:(x^{\prime },p_{2}^{\prime
},q_{2}^{\prime })\longrightarrow (X^{\prime },P_{2}^{\prime },Q_{2}^{\prime
})$ and
\begin{eqnarray*}
(\widehat{x},\widehat{p},\widehat{q}) &=&(x-x^{\prime },p_{2}-p_{2}^{\prime
},q_{2}-q_{2}^{\prime }), \\
(\widehat{X},\widehat{P},\widehat{Q}) &=&(X-X^{\prime },P_{2}-P_{2}^{\prime
},Q_{2}-Q_{2}^{\prime })
\end{eqnarray*}%
\begin{equation*}
\left\{
\begin{array}{rcl}
\mathrm{d}\widehat{X} & = & \left[ \alpha _{0}\widehat{\mathbb{B}}\left(
\widehat{X},\widehat{P},\widehat{Q}\right) +\delta \widehat{\mathbb{B}}%
\left( \widehat{x},\widehat{p},\widehat{q}\right) \right] \mathrm{d}t \\
&  & +\left[ \alpha _{0}\widehat{\mathbb{C}}\left( \widehat{X},\widehat{P},%
\widehat{Q}\right) +\delta \widehat{\mathbb{C}}\left( \widehat{x},\widehat{p}%
,\widehat{q}\right) \right] \mathrm{d}W\left( t\right) , \\
-\mathrm{d}\widehat{P} & = & \left[ \alpha _{0}\widehat{\mathbb{F}}\left(
\widehat{X},\widehat{P},\widehat{Q}\right) +\delta \widehat{\mathbb{F}}%
\left( \widehat{x},\widehat{p},\widehat{q}\right) \right] \mathrm{d}t-%
\widehat{Q}\mathrm{d}W\left( t\right) \\
\widehat{X}\left( 0\right) & = & 0,\text{ }\widehat{P}\left( T\right)
=-\alpha _{0}\Phi _{2}\widehat{X}\left( T\right) -\delta \Phi _{2}\widehat{x}%
\left( T\right) ,%
\end{array}%
\right.
\end{equation*}%
with
\begin{equation*}
\left\{
\begin{array}{l}
\widehat{\mathbb{B}}\left( \widehat{X},\widehat{P},\widehat{Q}\right) =%
\mathbb{B}\left( X,P_{2},Q_{2}\right) -\mathbb{B}\left( X,P_{2}^{\prime
},Q_{2}^{\prime }\right) , \\
\widehat{\mathbb{C}}\left( \widehat{X},\widehat{P},\widehat{Q}\right) =%
\mathbb{C}\left( X,P_{2},Q_{2}\right) -\mathbb{C}\left( X,P_{2}^{\prime
},Q_{2}^{\prime }\right) , \\
\widehat{\mathbb{F}}\left( \widehat{X},\widehat{P},\widehat{Q}\right) =%
\mathbb{F}\left( X,P_{2},Q_{2}\right) -\mathbb{F}\left( X,P_{2}^{\prime
},Q_{2}^{\prime }\right) ,%
\end{array}%
\right.
\end{equation*}%
Applying It\^{o} formula to $\big<\widehat{P},\widehat{X}\big>$ and taking
expectations on both sides:%
\begin{eqnarray*}
&&\mathbb{E}\left\langle \widehat{X}\left( T\right) ,-\alpha _{0}\Phi _{2}%
\widehat{X}\left( T\right) -\delta \Phi _{2}\widehat{x}\left( T\right)
\right\rangle \\
&=&\mathbb{E}\bigg \{\int_{0}^{T}\bigg [\left\langle \widehat{X}\left(
s\right) ,-\alpha _{0}\widehat{\mathbb{F}}\left( \widehat{X}\left( s\right) ,%
\widehat{P}\left( s\right) ,\widehat{Q}\left( s\right) \right) \right\rangle
+\left\langle \widehat{X}\left( s\right) ,-\delta \widehat{\mathbb{F}}\left(
\widehat{x}\left( s\right) ,\widehat{p}\left( s\right) ,\widehat{q}\left(
s\right) \right) \right\rangle \\
&&+\left\langle \widehat{P}\left( s\right) ,\alpha _{0}\widehat{\mathbb{B}}%
\left( \widehat{X}\left( s\right) ,\widehat{P}\left( s\right) ,\widehat{Q}%
\left( s\right) \right) \right\rangle +\left\langle \widehat{P}\left(
s\right) ,\delta \widehat{\mathbb{B}}\left( \widehat{x}\left( s\right) ,%
\widehat{p}\left( s\right) ,\widehat{q}\left( s\right) \right) \right\rangle
\\
&&+\left\langle \widehat{Q}\left( s\right) ,\alpha _{0}\widehat{\mathbb{C}}%
\left( \widehat{X}\left( s\right) ,\widehat{P}\left( s\right) ,\widehat{Q}%
\left( s\right) \right) \right\rangle +\left\langle \widehat{Q}\left(
s\right) ,\delta \widehat{\mathbb{C}}\left( \widehat{x}\left( s\right) ,%
\widehat{p}\left( s\right) ,\widehat{q}\left( s\right) \right) \right\rangle %
\bigg ]\mathrm{d}s\bigg \}
\end{eqnarray*}
Rearranging the above terms, we have%
\begin{eqnarray*}
&&\alpha _{0}\mathbb{E}\left\langle \widehat{X}\left( T\right) ,\Phi _{2}%
\widehat{X}\left( T\right) \right\rangle +\alpha _{0}\mathbb{E}\bigg \{%
\int_{0}^{T}\bigg [\left\langle \widehat{X}\left( s\right) ,-\widehat{%
\mathbb{F}}\left( \widehat{X}\left( s\right) ,\widehat{P}\left( s\right) ,%
\widehat{Q}\left( s\right) \right) \right\rangle \\
&&+\left\langle \widehat{P}\left( s\right) ,\widehat{\mathbb{B}}\left(
\widehat{X}\left( s\right) ,\widehat{P}\left( s\right) ,\widehat{Q}\left(
s\right) \right) \right\rangle +\left\langle \widehat{Q}\left( s\right) ,%
\widehat{\mathbb{C}}\left( \widehat{X}\left( s\right) ,\widehat{P}\left(
s\right) ,\widehat{Q}\left( s\right) \right) \right\rangle \bigg ]\mathrm{d}s%
\bigg \} \\
&=&\delta \mathbb{E}\bigg \{\int_{0}^{T}\bigg [\left\langle \widehat{X}%
\left( s\right) ,\widehat{\mathbb{F}}\left( \widehat{x}\left( s\right) ,%
\widehat{p}\left( s\right) ,\widehat{q}\left( s\right) \right) \right\rangle
-\left\langle \widehat{P}\left( s\right) ,\widehat{\mathbb{B}}\left(
\widehat{x}\left( s\right) ,\widehat{p}\left( s\right) ,\widehat{q}\left(
s\right) \right) \right\rangle \\
&&-\left\langle \widehat{Q}\left( s\right) ,\delta \widehat{\mathbb{C}}%
\left( \widehat{x}\left( s\right) ,\widehat{p}\left( s\right) ,\widehat{q}%
\left( s\right) \right) \right\rangle \bigg ]\mathrm{d}s\bigg \}-\delta
\mathbb{E}\left\langle \widehat{X}\left( T\right) ,\Phi _{2}\widehat{x}%
\left( T\right) \right\rangle
\end{eqnarray*}

Hence,%
\begin{eqnarray}
&&\alpha _{0}\mathbb{E}\left\vert \Phi _{2}^{\frac{1}{2}}\widehat{X}\left(
T\right) \right\vert ^{2}+\alpha _{0}\mathbb{E}\bigg \{\int_{0}^{T}\left%
\vert Q_{2}^{\frac{1}{2}}\left( s\right) \widehat{X}\left( s\right)
\right\vert ^{2}\mathrm{d}s\bigg \}  \notag \\
&\leq &\alpha _{0}\mathbb{E}\left\langle \widehat{X}\left( T\right) ,\Phi
_{2}\widehat{X}\left( T\right) \right\rangle +\alpha _{0}\mathbb{E}\bigg \{%
\int_{0}^{T}\bigg [\left\langle \widehat{X}\left( s\right) ,-\widehat{%
\mathbb{F}}\left( \widehat{X}\left( s\right) ,\widehat{P}\left( s\right) ,%
\widehat{Q}\left( s\right) \right) \right\rangle  \notag \\
&&+\left\langle \widehat{P}\left( s\right) ,\widehat{\mathbb{B}}\left(
\widehat{X}\left( s\right) ,\widehat{P}\left( s\right) ,\widehat{Q}\left(
s\right) \right) \right\rangle +\left\langle \widehat{Q}\left( s\right) ,%
\widehat{\mathbb{C}}\left( \widehat{X}\left( s\right) ,\widehat{P}\left(
s\right) ,\widehat{Q}\left( s\right) \right) \right\rangle \bigg ]\mathrm{d}s%
\bigg \}  \notag \\
&=&\delta \mathbb{E}\bigg \{\int_{0}^{T}\bigg [\left\langle \widehat{X}%
\left( s\right) ,\widehat{\mathbb{F}}\left( \widehat{x}\left( s\right) ,%
\widehat{p}\left( s\right) ,\widehat{q}\left( s\right) \right) \right\rangle
-\left\langle \widehat{P}\left( s\right) ,\widehat{\mathbb{B}}\left(
\widehat{x}\left( s\right) ,\widehat{p}\left( s\right) ,\widehat{q}\left(
s\right) \right) \right\rangle  \notag \\
&&-\left\langle \widehat{Q}\left( s\right) ,\delta \widehat{\mathbb{C}}%
\left( \widehat{x}\left( s\right) ,\widehat{p}\left( s\right) ,\widehat{q}%
\left( s\right) \right) \right\rangle \bigg ]\mathrm{d}s\bigg \}-\delta
\mathbb{E}\left\langle \widehat{X}\left( T\right) ,\Phi _{2}\widehat{x}%
\left( T\right) \right\rangle  \notag \\
&\leq &\delta \mathbb{E}\bigg \{\int_{0}^{T}\bigg [\left\vert \widehat{X}%
\left( s\right) \right\vert ^{2}+\left\vert \widehat{P}\left( s\right)
\right\vert ^{2}+\left\vert \widehat{Q}\left( s\right) \right\vert ^{2}\bigg
]\mathrm{d}s\bigg \}  \notag \\
&&+\delta \mathbb{E}\bigg \{\int_{0}^{T}\bigg [\left\vert \widehat{x}\left(
s\right) \right\vert ^{2}+\left\vert \widehat{p}\left( s\right) \right\vert
^{2}+\left\vert \widehat{q}\left( s\right) \right\vert ^{2}\bigg ]\mathrm{d}s%
\bigg \}  \notag \\
&&+\delta C\left( \left\vert \widehat{x}\left( T\right) \right\vert
^{2}+\left\vert \widehat{X}\left( T\right) \right\vert ^{2}\right) .
\label{ineq1}
\end{eqnarray}%
We point out that, the first inequality uses the monotonicity property of $%
\varphi (p,q)$ (Proposition \ref{mono}). The second inequality is due to the
basic geometric inequality and Lipschitz property of projection operator
(Proposition \ref{Plip}).

Then, by standard estimates of BSDE:
\begin{eqnarray}
&&\mathbb{E}\bigg \{\int_{0}^{T}\bigg [\left\vert \widehat{P}\left( s\right)
\right\vert ^{2}+\left\vert \widehat{Q}\left( s\right) \right\vert ^{2}\bigg
]\mathrm{d}s\bigg \}  \notag \\
&\leq &\delta C\mathbb{E}\bigg \{\int_{0}^{T}\bigg [\left\vert \widehat{x}%
\left( s\right) \right\vert ^{2}+\left\vert \widehat{p}\left( s\right)
\right\vert ^{2}+\left\vert \widehat{q}\left( s\right) \right\vert ^{2}\bigg
]\mathrm{d}s\bigg \}+\delta C\mathbb{E}\left\vert \widehat{x}\left( T\right)
\right\vert ^{2}  \notag \\
&&+C\alpha _{0}\left( \mathbb{E}\left\vert \Phi _{2}^{\frac{1}{2}}\widehat{X}%
\left( T\right) \right\vert ^{2}+\mathbb{E}\bigg \{\int_{0}^{T}\left\vert
Q_{2}^{\frac{1}{2}}\left( s\right) \widehat{X}\left( s\right) \right\vert
^{2}\mathrm{d}s\bigg \}\right) .  \label{ineq2}
\end{eqnarray}

Next, by the standard estimate of forward SDEs:%
\begin{eqnarray}
&&\mathbb{E}\bigg \{\int_{0}^{T}\left\vert \widehat{X}\left( s\right)
\right\vert ^{2}\mathrm{d}s\bigg \}+\mathbb{E}\left\vert \widehat{X}\left(
T\right) \right\vert ^{2}  \notag \\
&\leq &\delta C\mathbb{E}\bigg \{\int_{0}^{T}\bigg [\left\vert \widehat{x}%
\left( s\right) \right\vert ^{2}+\left\vert \widehat{p}\left( s\right)
\right\vert ^{2}+\left\vert \widehat{q}\left( s\right) \right\vert ^{2}\bigg
]\mathrm{d}s\bigg \}  \notag \\
&&+C\mathbb{E}\left[ \int_{0}^{T}\bigg [\left\vert \widehat{P}\left(
s\right) \right\vert ^{2}+\left\vert \widehat{Q}\left( s\right) \right\vert
^{2}\bigg ]\mathrm{d}s\right] +\delta C\mathbb{E}\left\vert \widehat{x}%
\left( T\right) \right\vert ^{2}.  \label{ineq3}
\end{eqnarray}%
Based on the above estimates (\ref{ineq1})-(\ref{ineq3}), we see the mapping
$I$ satisfying%
\begin{equation*}
\mathbb{E}\int_{0}^{T}\left( |\widehat{X}_{s}|^{2}+|\widehat{P}_{s}|^{2}+|%
\widehat{Q}_{s}|^{2}\right) ds+\mathbb{E}|\widehat{X}_{T}|^{2}\leq K\delta
\left( \mathbb{E}\int_{0}^{T}\left( |\widehat{x}_{s}|^{2}+|\widehat{p}%
_{s}|^{2}+|\widehat{q}_{s}|^{2}\right) ds+\mathbb{E}|\widehat{x}%
_{T}|^{2}\right) .
\end{equation*}%
It follows the mapping is a contraction and the existence follows
immediately using the arguments presented in \cite{hp} and \cite{PW1999}.
~\hfill $\Box $

\section{Discussion on Riccati equation}

Recall the stochastic Hamilton system (taken from Tang \cite{Tang03}) is
given by
\begin{equation*}
\left\{
\begin{array}{rcl}
\mathrm{d}x\left( t\right) & = & \left( A\left( t\right) x\left( t\right)
+B\left( t\right) u\left( t\right) \right) \mathrm{d}t+\sum_{i=1}^{d}\left(
C^{i}\left( t\right) x\left( t\right) +D^{i}\left( t\right) u\left( t\right)
\right) \mathrm{d}W^{i}\left( t\right) , \\
u\left( t\right) & = & -N^{-1}\left( t\right) \left[ B^{\top }\left(
t\right) y\left( t\right) +\sum_{i=1}^{d}D^{i}\left( t\right) ^{\top
}z^{i}\left( t\right) \right] , \\
-\mathrm{d}y\left( t\right) & = & \left[ A^{\top }\left( t\right) y\left(
t\right) +\sum_{i=1}^{d}C^{i}\left( t\right) ^{\top }z^{i}\left( t\right)
+Q\left( t\right) x\left( t\right) \right] \mathrm{d}t-\sum_{i=1}^{d}z^{i}%
\left( t\right) \mathrm{d}W^{i}\left( t\right) \\
x\left( \tau \right) & = & h\in L^{2}\left( \Omega ,\mathcal{F}_{\tau },P;%
\mathbb{R}^{n}\right) ,\text{ }y\left( T\right) =Mx\left( T\right) .%
\end{array}%
\right.
\end{equation*}%
Inserting $u\left( \cdot \right) $ into the first equation, we have
\begin{equation*}
\left\{
\begin{array}{rcl}
\mathrm{d}x & = & \left( Ax-BN^{-1}\left[ B^{\top }y+\sum_{i=1}^{d}\left(
D^{i}\right) ^{\top }z^{i}\right] \right) \mathrm{d}t \\
&  & +\sum_{i=1}^{d}\left( C^{i}x-D^{i}N^{-1}\left[ B^{\top
}y+\sum_{i=1}^{d}\left( D^{i}\right) ^{\top }z^{i}\right] \right) \mathrm{d}%
W^{i}\left( t\right) , \\
-\mathrm{d}y & = & \left[ A^{\top }y+\sum_{i=1}^{d}\left( C^{i}\right)
^{\top }z^{i}+Qx\right] \mathrm{d}t-\sum_{i=1}^{d}z^{i}\mathrm{d}W^{i}\left(
t\right) \\
x\left( \tau \right) & = & h\in L^{2}\left( \Omega ,\mathcal{F}_{\tau },P;%
\mathbb{R}^{n}\right) ,\text{ }y\left( T\right) =Mx\left( T\right) .%
\end{array}%
\right.
\end{equation*}%
A formal approach to derive the associated Riccati equation from the
stochastic Hamilton system a priori assumes that there is a semi-martingale $%
K $ of the form%
\begin{equation*}
K\left( t\right) =K\left( 0\right) -\int_{0}^{t}K_{1}\left( s\right) \mathrm{%
d}s+\int_{0}^{t}\sum_{i=1}^{d}L^{i}\left( s\right) \mathrm{d}W\left(
s\right) ,\text{ }0\leq t\leq T.
\end{equation*}%
such that
\begin{equation*}
y\left( t\right) =K\left( t\right) x\left( t\right) .
\end{equation*}%
Then, applying It\^{o}'s formula to $K\left( t\right) x\left( t\right) $, we
have%
\begin{eqnarray}
&&K\left( Ax-BN^{-1}\left[ B^{\top }y+\sum_{i=1}^{d}\left( D^{i}\right)
^{\top }z^{i}\right] \right) \mathrm{d}t  \notag \\
&&+K\sum_{i=1}^{d}\left( C^{i}x-D^{i}N^{-1}\left[ B^{\top
}y+\sum_{i=1}^{d}\left( D^{i}\right) ^{\top }z^{i}\right] \right) \mathrm{d}%
W^{i}\left( t\right)  \notag \\
&&-K_{1}x\mathrm{d}t+\sum_{i=1}^{d}L^{i}x\mathrm{d}W^{i}\left( t\right)
\notag \\
&&+\sum_{i=1}^{d}L^{i}\left( s\right) \sum_{i=1}^{d}\left(
C^{i}x-D^{i}N^{-1} \left[ B^{\top }y+\sum_{i=1}^{d}\left( D^{i}\right)
^{\top }z^{i}\right] \right) \mathrm{d}t  \notag \\
&=&\mathrm{d}y\left( t\right)  \notag \\
&=&-\left[ A^{\top }y+\sum_{i=1}^{d}\left( C^{i}\right) ^{\top }z^{i}+Qx%
\right] \mathrm{d}t+\sum_{i=1}^{d}z^{i}\mathrm{d}W^{i}\left( t\right) .
\label{tang}
\end{eqnarray}%
It follows that
\begin{equation*}
z^{i}=K\left( C^{i}x-D^{i}N^{-1}B^{\top }Kx-D^{i}N^{-1}\left( D^{i}\right)
^{\top }z^{i}\right) +L^{i}x,\text{ }1\leq i\leq d.
\end{equation*}%
Immediately,
\begin{equation}
z^{i}=Z^{i}x,  \label{zz}
\end{equation}%
where
\begin{equation*}
Z^{i}=\left( I+KD^{i}N^{-1}\left( D^{i}\right) ^{\top }\right) ^{-1}\left(
KC^{i}-KD^{i}N^{-1}B^{\top }K+L^{i}\right) .
\end{equation*}%
Substituting (\ref{zz})\footnote{%
In order to get the standard form (3.1) in Tang \cite{Tang03}, the rest
proceeding will employ the well-known matrix inverse formula, that is,
\begin{equation}
\left( A+BCD\right) ^{-1}=A^{-1}-A^{-1}B\left( DA^{-1}B+C^{-1}\right)
^{-1}DA^{-1}
\end{equation}%
where $A^{-1},$ $C^{-1},$ and either $\left( A+BCD\right) ^{-1}$ or $\left(
DA^{-1}B+C^{-1}\right) ^{-1}$ are assumed to exist.} into (\ref{tang}), we
identify the integrands of the Lebesgue integral. As a consequence%
\begin{eqnarray}
&&KA-K\underset{\boldsymbol{B}_{1}}{\underbrace{BN^{-1}B^{\top }}}K-K%
\underset{\boldsymbol{B}_{2}}{\underbrace{BN^{-1}\sum_{i=1}^{d}\left(
D^{i}\right) ^{\top }}}Z^{i}-K_{1}  \notag \\
&&+\sum_{i=1}^{d}L^{i}\left( s\right) \sum_{i=1}^{d}\left( C^{i}-\underset{%
\boldsymbol{D}_{1}}{\underbrace{D^{i}N^{-1}B^{\top }}}K-\underset{%
\boldsymbol{D}_{2}}{\underbrace{D^{i}N^{-1}\left( D^{i}\right) ^{\top }}}%
Z^{i}\right)  \notag \\
&=&-A^{\top }K-\sum_{i=1}^{d}\left( C^{i}\right) ^{\top }Z^{i}-Q.
\label{tang2}
\end{eqnarray}%
Now
\begin{eqnarray*}
K_{1} &=&A^{\top }K+\sum_{i=1}^{d}\left( C^{i}\right) ^{\top }Z^{i}+Q+KA-K%
\boldsymbol{B}_{1}K-K\boldsymbol{B}_{2}Z^{i} \\
&&+\sum_{i=1}^{d}L^{i}\left( s\right) \sum_{i=1}^{d}\left( C^{i}-\boldsymbol{%
D}_{1}K-\boldsymbol{D}_{2}Z^{i}\right) .
\end{eqnarray*}%
This is another expression for Riccati equation (3.1) in Tang \cite{Tang03}.

\end{document}